\documentclass{amsart}

\usepackage{amsthm}
\usepackage{psfrag}
\usepackage{graphicx}
\usepackage{amssymb}
\newtheorem{teo}{Theorem}[section]
\newtheorem{lemma}[teo]{Lemma}
\newtheorem{prop}[teo]{Proposition}
\newtheorem{cor}[teo]{Corollary}

\theoremstyle{definition}
\newtheorem{defi}[teo]{Definition}
\newtheorem{rem}[teo]{Remark}
\newtheorem{example}[teo]{Example}

\theoremstyle{remark}
\newtheorem{prof}[teo]{Proof of}

\newcommand{\qbin}[2]{\left[\begin{array}{c}
#1 \\
#2 \end{array}\right]}

\def\mc{\mathbb{C}}
\def\mr{\mathbb{R}}
\def\mz{\mathbb{Z}}

\def\nns{\negthickspace}

\newcommand{\smalltheta}[5]{\raisebox{-#5cm}{\psfrag{a}{\footnotesize{$#1$}}\psfrag{b}{\footnotesize{$#2$}}\psfrag{c}{\footnotesize{$#3$}}\includegraphics[width=#4cm]{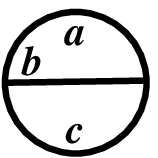}}}
\newcommand{\smallunknot}[3]{\raisebox{-#3cm}{\psfrag{a}{\footnotesize{$#1$}}\includegraphics[width=#2cm]{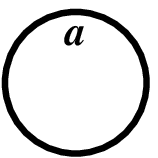}}}
\newcommand{\smalltet}[8]{\raisebox{-#8cm}{\psfrag{a}{\footnotesize{$#1$}}\psfrag{b}{\footnotesize{$#2$}}\psfrag{c}{\footnotesize{$#3$}}\psfrag{d}{\footnotesize{$#4$}}\psfrag{e}{\footnotesize{$#5$}}\psfrag{f}{\footnotesize{$#6$}}\includegraphics[width=#7cm]{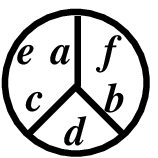}}}
\newcommand{\smalltetcross}[8]{\raisebox{-#8cm}{\psfrag{a}{\footnotesize{$#1$}}\psfrag{b}{\footnotesize{$#2$}}\psfrag{c}{\footnotesize{$#3$}}\psfrag{d}{\footnotesize{$#4$}}\psfrag{e}{\footnotesize{$#5$}}\psfrag{f}{\footnotesize{$#6$}}\includegraphics[width=#7cm]{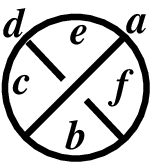}}}
\newcommand{\JW}[1]{\raisebox{-.6cm}{\psfrag{JW}{\footnotesize{$#1$}}\includegraphics[width=.9cm]{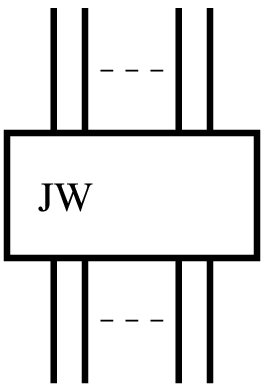}}}
\newcommand{\sigmahat}[1]{\raisebox{-.6cm}{\psfrag{sigma}{\footnotesize{$#1$}}\includegraphics[width=.9cm]{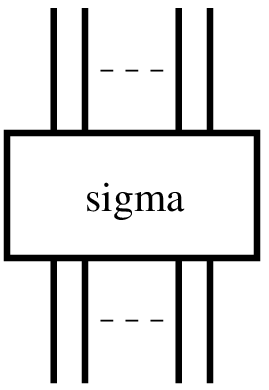}}}
\newcommand{\poscross}[2]{\raisebox{-.4cm}{\psfrag{a}{\footnotesize{$#1$}}\psfrag{b}{\footnotesize{$#2$}}\includegraphics[width=.9cm]{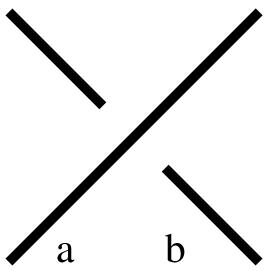}}}

\newcommand{\dcupcap}[1]{\raisebox{-.4cm}{\psfrag{a}{\footnotesize{$#1$}}\includegraphics[width=.9cm]{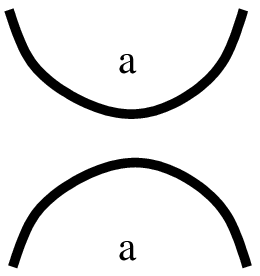}}}
\newcommand{\didentity}[1]{\raisebox{-.5cm}{\psfrag{a}{\footnotesize{$#1$}}\includegraphics[width=.6cm]{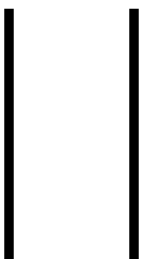}}}

\begin{document}
\title{Integrality of Kauffman brackets of trivalent graphs}

\author[Costantino]{Francesco Costantino}
\address{Institut de Recherche Math\'ematique Avanc\'ee\\
  Rue Ren\'e Descartes 7\\
  67084 Strasbourg, France}
\email{costanti@math.u-strasbg.fr}

\begin{abstract}
We show that Kauffman brackets of colored framed graphs (also known as quantum spin networks) can be renormalized to a Laurent polynomial with integer coefficients by multiplying it by a coefficient which is a product of quantum factorials depending only on the abstract combinatorial structure of the graph. Then we compare the shadow-state sums and the state-sums based on $R$-matrices and Clebsch-Gordan symbols, reprove their equivalence and comment on the integrality of the weight of the states. We also provide short proofs of most of the standard identities satisfied by quantum $6j$-symbols of $U_q(sl_2)$.
\end{abstract}

\maketitle

\tableofcontents
\section{Introduction}
The family of $U_q(sl_2)$-quantum invariants of knotted objects in $S^3$ as knots, links and more in general trivalent graphs, can be defined via the recoupling theory (\cite{KL}) as well as via the theory of representations of the quantum group $U_q(sl_2)$ (\cite{KM},\cite{RT}) and the so-called theory of ``shadows" (\cite{Tu}, Chapter IX). These invariants are defined for framed objects (framed links or graphs, see Definition \ref{def:ktg}) equipped with a ``coloring" on the set of $1$-dimensional strata of the object satisfing certain admissibility conditions (Definition \ref{def:admcol}), and take values in $\mathbb{Q}(q^{\frac{1}{2}})$. As customary in the literature, we shall call these invariants the ``Kauffman brackets" of the colored, framed graph $G$ and denote them by $\langle   G,col\rangle$. 
In particular, if $G$ is a framed knot $k$ the set of admissible colorings is the set of half-natural numbers (here we use the so-called ``spin" notation) and they coincide with the unreduced colored Jones polynomials of the knot:  $\langle   k,s\rangle =J_{2s+1}(k)$. 

Although the definition of Kauffman brackets via recoupling theory is simple and appealing, the definition based on the theory of representations of $U_q(sl_2)$ happens to be more useful for our purposes. In Section \ref{sec:kauffman} we will sketch the proof of the equivalence of the two definitions (the relations have been already studied by S. Piunikhin \cite{Pi} and is fully detailed in \cite{CS}). 

It is known that, in general, $\langle   G,col\rangle $ is a rational function of the variable $q^{\frac{1}{2}}$ (there are various notations in the literature, for instance our $q^{\frac{1}{2}}$ is $A$ in \cite{KL}).
If $L$ is a framed link, it was shown by T. Le (\cite{Le}) that, up to a factor of the form $ q^{\pm\frac{n}{2}}$, $\langle  L,n\rangle $ is a Laurent polynomial in $q$ (actually in \cite{Le} a much stronger result is proved which holds for general polynomial invariants issued from quantum group representations).
 
On contrast it is well known that $\langle  G,col\rangle $ is not in general a Laurent polynomial if $G$ is a trivalent graph. 
The main result of the present paper is Theorem \ref{teo:main}, restated here in a simpler form (we refer to Section \ref{sec:main} for the notation):\begin{teo}[Integrality of the renormalized Kauffman brackets]\label{teo:intro}
There exist $m,n\in \mz$ such that:
$$
\langle  \langle  G,col\rangle \rangle\doteqdot \langle  G,col\rangle  \frac{\prod [2col(e)]!}{\prod[a_v+b_v-c_v]![b_v+c_v-a_v]![c_v+a_v-b_v]!}
  \in (\sqrt{-1})^mq^{\frac{n}{4}}\mathbb{Z}[q,q^{-1}]
$$
where the products are taken over the non-closed edges $e$ of $G$ and the vertices $v$ of $G$.  
\end{teo}

 It turns out that $\langle  \langle  G,col\rangle \rangle =\langle  G,col\rangle $ if $G$ is a link.
This normalization was proposed and conjectured to be integral by S. Garoufalidis and R. Van der Veen (in \cite{GV}, where they also proved the integrality in the classical case when $q=\pm 1$) in order to define generating function for classical spin networks evaluations. We hope that this result will allow further development in that direction and in the understanding of the categorification of $U_q(sl_2)$-quantum invariants for general knotted objects.

The last sections are almost independent from the preceding ones. In Section \ref{sec:shadows} we recall the definition of shadow-state sums to compute $\langle\langle G,col\rangle\rangle$, give a new self-contained proof of the equivalence (first proved in \cite{KR}) between the shadow-state formulation and the $R$-matrix formulation of the invariants, and comment on the non-integrality of the single shadow-state weights.
Using shadow state-sums we also provide short proofs of the most famous identities for $6j$-symbols (e.g. Racah, Biedenharn-Elliot, orthogonality). 
In the last section we will quickly comment on the case when $G$ has non-empty boundary and on the algebraic meaning of the shadow-state sums with respect to the state-sum based on $R$-matrices and Clebsch-Gordan symbols.

\subsection{Structure of the paper} 
In Section \ref{sec:kauffman} we will recall the definition of Kauffman brackets of colored graphs and  the basic facts on representation theory of $U_q(sl_2)$ (for generic $q$). We then show how to compute Kauffman brackest via morphisms associated to tangles, and provide explicit formulas for the elementary morphisms. In Section \ref{sec:main} we will define $\langle  G,col\rangle $ and prove its integrality. Section \ref{sec:shadows} is almost independent from the first sections (basically it depends only on Lemma \ref{lem:properties}); there we explain how to compute $\langle\langle G,col\rangle\rangle$ via shadow state-sums and provide short proofs of some well known identities for $6j$-symbols. In Section \ref{sec:Rvs6j} we comment on the algebraic meaning of shadow-state-sums. 
\subsection{Acnowledgements}
I wish to thank Fran\c cois Gueritaud, Roland Van der Veen and Vladimir Turaev for the comments and suggestions they gave me. 
This work was supported by the French ANR project ANR-08-JCJC-0114-01.
\section{Kauffman brackets via representations of $U_q(sl_2)$}\label{sec:kauffman}
\subsection{The definition of Kauffman brackets}\label{sub:kauffman}
\begin{defi}[KTG]\label{def:ktg}
A $Knotted\ Trivalent\ Graph$ (KTG) is a finite trivalent graph $G\subset S^3$ equipped with a ``framing", i.e. the germ of an orientable smooth surface $S\subset S^3$ such that $S$ retracts on $G$. 
\end{defi}
\begin{rem}
Note that this is not a ``fat graph" as $S$ is required to exist around all $G$ and not only around its vertices. On contrast we require $G$ to be embedded in $S^3$. 
Also, starting from the end of Subsection \ref{sub:kauffmanUq} we will drop the assumption for $S$ to be orientable. \end{rem}
In order to specify a framing $S$ on a graph $G$ we will only specify (via thin lines as in the leftmost drawing of Figure \ref{fig:correspondence}) the edges around which it twists with respect to the blackboard framing in a diagram of $G$, implicitly assuming that $S$ will be lying horizontally (i.e. parallel to the blackboard) around $G$ out of these twists.
Let us also remark that if $D$ is a diagram of $G$ there is a framing $S_D$ (called the \emph{blackboard framing}) induced on $G$ simply by considering a surface containing $G$ and lying almost parallel to the projection plane. Pulling back the orientation of $\mr^2$ shows that $S_D$ is orientable. The following is a converse: 
\begin{lemma}\label{lem:orientableframing}
If $G$ is a framed graph and $S$ is and orientable framing on $G$ then there exists a diagram $D$ of $G$ such that $S_D=S$.
\end{lemma}
\begin{prf}{1}{The idea of the proof is to fix a diagram $D$ and count the number of half twists of difference on each edge of $G$ between $S_D$ and $S$. The reduction mod $2$ of these numbers  forms an explicit cochain in $H^1(G;\mz_2)$ which is null cohomologous because $S$ and $S_D$ are orientable. The coboundary reducing it to the $0$ cochain corresponds to a finite number of moves as those in Lemma \ref{lem:halftwist} which change $D$ and isotope $G$ into a position such that the number of half twists of difference between $S$ and $S_D$ is even on every edge. Then up to adding a suitable number of kinks to each edge of $G$ this difference can be reduced to $0$ everywhere.
}\end{prf}

Let now $G$ be a KTG, $E$ the set of its edges, $V$ the set of its vertices: 
\begin{defi}[Admissible coloring]\label{def:admcol}
An \emph{admissible coloring} of $G$ is a map $col:E\to \frac{\mathbb{N}}{2}$ (whose values are called \emph{colors}) such that $\forall v\in V$ the following conditions are satisfied:
\begin{enumerate}
\item $a_v+b_v+c_v\in \mathbb{N}$
\item $a_v+b_v\geq c_v,\ b_v+c_v\geq a_v,\ c_v+a_v\geq b_v$
\end{enumerate}
where $a_v,b_v,c_v$ are the colors of the edges touching $v$. 
\end{defi}

Let us now fix diagram $D$ of $G$ such that the blackboard framing coincides with that of $G$ (it exists by Lemma \ref{lem:orientableframing}), and an admissible coloring $col$ on $G$, and recall how the Kauffmann bracket $<G,col>$ is defined.

Let $q\in\mc$, $[n]\doteqdot\frac{q^n-q^{-n}}{q-q^{-1}}$ and $[n]!\doteqdot \prod_{j=1}^n [j],\ [0]!=1$. Let also $\qbin{n}{k}\doteqdot \frac{[n]!}{[k]![n-k]!}$.
If $G$ is a framed link $L$ and $col$ is $\frac{1}{2}$ on all the components of $G$, the Kauffman bracket $<L,\frac{1}{2}>\in \mathbb{Z}[q^{\pm\frac{1}{2}}]$, is defined by applying recursively Kauffman's rules to $D$:
$$\poscross{}{}=q^{\frac{1}{2}}\didentity{}+q^{-\frac{1}{2}}\dcupcap{} \ \ \ \ \  \rm{and}\ \  \ \ \smallunknot{}{.9}{.4}=-[2]
$$
To define the general $<G,col>$ for trivalent colored graphs let's first define the Jones-Wenzl projectors $JW_{2a}\in C(q^{\frac{1}{2}})[B(2a)]$:
$$JW_{2a}=\JW{\ \ 2a}\doteqdot\sum_{\sigma \in \mathfrak{S}_{2a}} \frac{q^{-a(2a-1)+\frac{3}{2}T(\sigma)}}{[2a]!}\sigmahat{\ \hat{\sigma}}$$
where $\hat{\sigma}$ is the positive braid containing the minimal number ($T(\sigma)$) of crossings and inducing the permutation $\sigma$ on its endpoints (it is a standard fact that such braid is well defined).
Actually $JW_{2a}$ is defined as an element of the Temperly-Lieb algebra, but for the purpose of this section we will just consider it as a formal sum of braids; in the next section a more precise interpretation will be provided.
One defines $<G,col>$ by the following algorithm:
\begin{enumerate}
\item Cable each edge $e$ of $G$ by $JW_{2col(e)}$, i.e., in $D$ replace an edge $e$ colored by $a$ by a formal sum of braids in $B(2a)$ according to the above definition of $JW_{2a}$:
$$\raisebox{-.5cm}{\psfrag{a}{\footnotesize{$a$}}\includegraphics[width=.3cm]{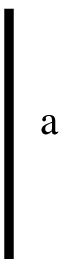}}\to \JW{\ 2a}$$

\item Around each vertex, connect the (yet free) endpoints of the so-obtained strands in the unique planar way without self returns:
$$
\raisebox{-.3cm}{\psfrag{a}{\footnotesize{$a$}}\psfrag{b}{\footnotesize{$b$}}\psfrag{c}{\footnotesize{$c$}}\includegraphics[width=1.5cm]{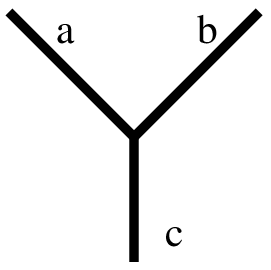}}\to \raisebox{-.5cm}{\psfrag{2a}{\footnotesize{$2a$}}\psfrag{2b}{\footnotesize{$2b$}}\psfrag{2c}{\footnotesize{$2c$}}\includegraphics[width=1.5cm]{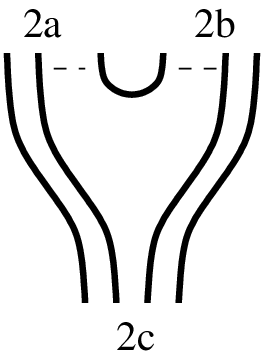}}$$
\item This way one associates to $(G,col)$ a formal sum with coefficients $c_i\in \mathbb{Q}(q^{1/2})$ of links $L_i$ contained in a small neighborhood of the framing of $G$ and therefore framed by annuli running parallel to it. Define $<G,col>\doteqdot \sum_i c_i<L_i,\frac{1}{2}>$.
\end{enumerate} 
\begin{teo}[Kauffman,\ \cite{KL}]
$<G,col>$ is an invariant up to isotopy of $(G,col)$.
\end{teo}

\subsection{Basic facts on $U_q(sl_2)$}
\begin{defi}
$U_q(sl_2)$ is the algebra generated by $E,F,K$ and $K^{-1}$ with relations:
$$[E,F]=\frac{K^2-K^{-2}}{q-q^{-1}},\  KE=qEK,\  KF=q^{-1}FK,\ KK^{-1}=K^{-1}K=1$$
Its Hopf algebra structure is given by: 
$$\Delta(E)=E\otimes K+K^{-1}\otimes E,\ \Delta(F)=F\otimes K+K^{-1}\otimes F,\ \Delta(K)=K\otimes K$$ 
$$S(E)=-qE,\ S(F)=-q^{-1}F,\ S(K)=K^{-1}$$ 
$$\epsilon(E)=\epsilon(F)=0,\ \epsilon(K)=1$$
\end{defi}
\begin{rem}
To make clear the relation with other notations note that ours is coherent with that of \cite{KM} after replacing their $s$ with $q$; our $q$ corresponds to $q=A^2$ in \cite{CS} and our $E,F$ respectively to $X$ and $Y$. 
\end{rem}
\begin{lemma}
For each $a\in \frac{\mathbb{N}}{2}$ there is a simple representation $V^a$ of $U_q(sl_2)$ of dimension $2a+1$ whose basis is $g^a_u,\ u=-a,-a+1,\ldots ,a$ and on which the action of $E,F,K$ is: 
$$E(g^a_u)=[a-u][a+u+1]g^a_{u+1},\ F(g^a_u)=g^a_{u-1},\ K(g^a_u)=q^u g^a_u$$
\end{lemma}
\begin{rem}
In \cite{CS} and \cite {KM} different bases $e^a_u$ and $f^{a}_u$ for $V^a$ where used. 
The changes of basis are in both cases diagonal and are: $f^a_m=\frac{[2a]!}{[a-u]!}g^a_u$ and $e^a_u=[a+u]!g^{a}_u$. 
\end{rem}
Let also recall that, by Clebsch-Gordan decomposition theorem, $V^a\otimes V^b$ is isomorphic to $V^{a+b}\oplus V^{a+b-1}\oplus \ldots \oplus V^{\vert a-b\vert }$.
Hence by Schur's lemma, the space $Hom(V^c,V^a\otimes V^b)$ has dimension $1$ if $(a,b,c)$ is admissible and $0$ otherwise; in the next subsection, for each three-uple $a,b,c$ we will choose explicit elements $Y_c^{a,b}\in Hom(V^c,V^a\otimes V^b)$ which are ``induced from the topology".
\subsection{Computing Kauffman brackets invariants via $U_q$}\label{sub:kauffmanUq}
\begin{figure}
\psfrag{a}{$a$}
\psfrag{b}{$b$}
\psfrag{c}{$c$}
\psfrag{j}{$a$}
\psfrag{cup}{$\cup_a$}
\psfrag{cap}{$\cap_a$}
\psfrag{R}{${}_b^aR$}
\psfrag{R-}{${}_b^aR_-$}
\psfrag{Y}{$Y^{a,b}_{c}$}
\psfrag{H}{$H_a$}
\psfrag{W}{$W^{a,b,c}$}
\psfrag{M}{$M^{a,b,c}$}
\psfrag{C}{$C^{a,b,c}$}
\psfrag{P}{$P^c_{a,b}$}
\includegraphics[width=12cm]{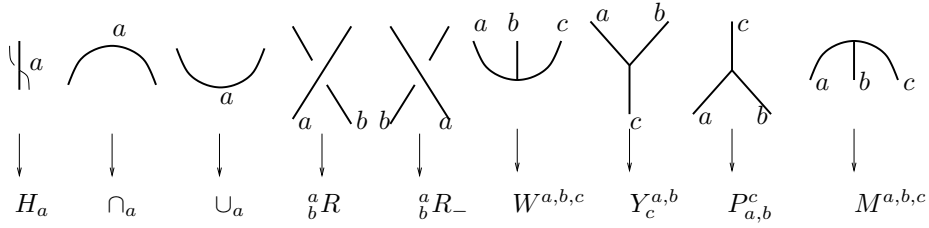}
\caption{The elementary graphs and the associated morphisms. In all the drawings except the leftmost, the framing is the blackboard framing.}\label{fig:correspondence}
\end{figure}

The standard construction of quantum invariants via the representation theory of $U_q$ (\cite {KM},\cite{RT}) allows one to associate to each diagram of an $(n,m)$-colored framed tangle $G$ (possibly containing some vertices) an operator between representations of $U_q$.
More explicitly if the bottom strands of $G$ are colored by $a_1,\ldots a_n$ and the top strands by $b_1,\ldots b_m$ then one can associate to a diagram $D$ of $G$ a morphisms $op(G,col,D):V_{a_1}\otimes\cdots \otimes V_{a_n}\to V_{b_1}\otimes\cdots \otimes V_{b_m}$. 
To do this, one defines the operators associated to each ``elementary" subdiagram (shown in Figure \ref{fig:correspondence}) equipped with any admissible coloring and then decomposes $D$ into a vertical stacking of these subdiagrams: $op(G,col,D)$ is then defined as a composition of the operators associated to the elementary blocks.
If one can choose the elementary operators so that the resulting morphism $op(G,col,D)$ does not depend on $D$, then in particular, if $G$ is closed, $op(G,col,D)$ will be an invariant up to isotopy of $(G,col)$ with values in $C(q^{\frac{1}{4}})$:
\begin{prop}\label{prop:kauffmanviauq}
There exist choices of operators for all admissible colorings of the elementary diagrams of Figure \ref{fig:correspondence} such that for each closed, colored KTG $(G,col)$, and each diagram $D$ of $G$ such that the framing of $G$ coincides with the blackboard framing, it holds $op(G,col,D)=<G,col>$.
\end{prop}
Let us define operators $\cup_{\frac{1}{2}}:V^0\to V^{\frac{1}{2}}\otimes V^{\frac{1}{2}}$ and $\cap_{\frac{1}{2}}:V^{\frac{1}{2}}\otimes V^{\frac{1}{2}}\to V^0$ which will allow us to compute $<G,col>$ whenever $G$ is a framed link and $col=\frac{1}{2}$.
We define it explicitly in the bases $g^{\frac{1}{2}}_u$ and $g^0_0$ by:
\begin{equation}\label{eq:cap}
\cap_{\frac{1}{2}}(g^{{\frac{1}{2}}}_{u}\otimes g^{{\frac{1}{2}}}_v)=\delta_{u,-v}\sqrt{-1}^{2u}q^u g^0_0
\end{equation}
\begin{equation}\label{eq:cup}
\cup_{\frac{1}{2}}(g^{0}_0)=\sum_{u=-{\frac{1}{2}}}^{{\frac{1}{2}}}\sqrt{-1}^{2u}q^u g^{\frac{1}{2}}_u\otimes g^{\frac{1}{2}}_{-u}
\end{equation}
Then we define the morphism ${}^{\frac{1}{2}}_{\frac{1}{2}}R:V^{\frac{1}{2}}\otimes V^{\frac{1}{2}}\to V^{\frac{1}{2}}\otimes V^{\frac{1}{2}}$ associated to a positive crossing of ${\frac{1}{2}}$-colored strands by $q^{{\frac{1}{2}}}Id+q^{-{\frac{1}{2}}}\cup_{\frac{1}{2}}\circ\cap_{\frac{1}{2}}$. Similarly, for a negative crossing we define ${}^{\frac{1}{2}}_{\frac{1}{2}}R_-:V^{\frac{1}{2}}\otimes V^{\frac{1}{2}}\to V^{\frac{1}{2}}\otimes V^{\frac{1}{2}}$ by $q^{-{\frac{1}{2}}}Id+q^{+{\frac{1}{2}}}\cup_{\frac{1}{2}}\circ\cap_{\frac{1}{2}}$.
Since $\cap_{\frac{1}{2}}\circ \cup_{\frac{1}{2}}=-[2]$, and by the definition of the crossing operator, it is evident that for each framed link $G$ colored by $\frac{1}{2}$ on each component and for any diagram $D$ of it, it will hold $op(G,col,D)=<G,col>$. Remark that $(\cap_{\frac{1}{2}}\otimes Id_{\frac{1}{2}})(Id_{\frac{1}{2}}\otimes \cup_{\frac{1}{2}})=Id_{\frac{1}{2}}$ and this encodes an isotopy invariance for planar diagrams; therefore for every framed $(n,m)$-tangle (implicitly colored by $\frac{1}{2}$) one has an associated morphism $(V^{\frac{1}{2}})^{\otimes n}\to (V^{\frac{1}{2}})^{\otimes m}$ well defined up to isotopy preserving the endpoints and the framing.

To treat the general case we define now $JW_{2a}\in End(V^{\otimes 2a}_{\frac{1}{2}})$ exactly as in Subsection \ref{sub:kauffman}, but now interpreting it as a morphism of representations of $U_q$ (via the definition of the $R$-operators we already gave).
The following holds (see \cite{CS}, Section 3.5):
\begin{teo}[Jones,Wenzl]
The operators $JW_{2a}$ are projectors over the unique submodule of $(V_{\frac{1}{2}})^{\otimes 2a}$ isomorphic to $V^a$.
\end{teo}
Therefore let us fix once and for all morphisms $\phi_a:V^a\to V^{\otimes 2a}_{\frac{1}{2}}$ by (see \cite{CS}, Definition 3.5.6, and recall that our base is related to that of \cite{CS} by a diagonal change $g^{a}_u=\frac{e^{a}_u}{[a+u]!}$): 
$$\phi_a(g^a_u)\doteqdot \frac{q^{\frac{a^2-u^2}{2}}}{[a+u]!}JW_{2a}((g^{\frac{1}{2}}_{\frac{1}{2}})^{\otimes (a+u)}\otimes (g^{\frac{1}{2}}_{-\frac{1}{2}})^{\otimes (a-u)})$$
and projectors $\mu_a:V^{\otimes 2a}_{\frac{1}{2}}\to V^a$, so that $\mu_a\circ \phi_a=Id_a$ and $\phi_a\circ \mu_a\circ JW_{2a}=JW_{2a}$. 
These operators allow to define the morphisms associated to a trivalent vertex $Y_c^{a,b}:V^c\to V^b\otimes V^c$ by $$\raisebox{-0.8cm}{\psfrag{a}{\footnotesize{$a$}}\psfrag{b}{\footnotesize{$b$}}\psfrag{c}{\footnotesize{$c$}}\includegraphics[width=1.4cm]{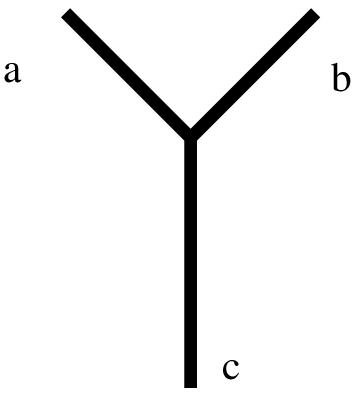}}\doteqdot Y_c^{a,b}\doteqdot(\mu_a\otimes \mu_b)\circ (JW_{2a}\otimes JW_{2b})\circ\raisebox{-1.0cm}{\psfrag{2a}{\footnotesize{$2a$}}\psfrag{2b}{\footnotesize{$2b$}}\psfrag{2c}{\footnotesize{$2c$}}\includegraphics[width=1.5cm]{opY.eps}}\circ JW_{2c}\circ\phi_c$$ where the drawing on the right represents the morphism $(V^{\frac{1}{2}})^{\otimes 2c}\to (V^{\frac{1}{2}})^{\otimes 2a}\otimes (V^{\frac{1}{2}})^{\otimes 2b}$ obtained by tensoring suitably $\cup_{\frac{1}{2}}$ and identity maps.
Similarly one can define maps $P_{a,b}^c:V^a\otimes V^b\to V^c$ by:
$$\raisebox{-0.7cm}{\psfrag{a}{\footnotesize{$a$}}\psfrag{b}{\footnotesize{$b$}}\psfrag{c}{\footnotesize{$c$}}\includegraphics[width=1.4cm]{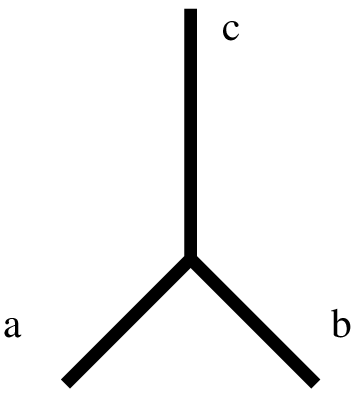}}\doteqdot P_{a,b}^c\doteqdot \mu_c\circ JW_{2c}\circ \raisebox{-1.0cm}{\psfrag{2a}{\footnotesize{$2a$}}\psfrag{2b}{\footnotesize{$2b$}}\psfrag{2c}{\footnotesize{$2c$}}\includegraphics[width=1.5cm]{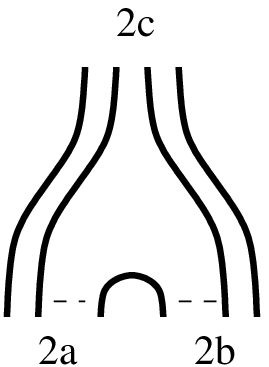}}\circ(JW_{2a}\otimes JW_{2b})\circ(\phi_a\otimes \phi_b)$$
From these, one gets new morphisms $\cup_a:V^0\to V^a\otimes V^a$ (as $Y_0^{(a,a)}$) and $\cap_{a}:V^a\otimes V^a\to V^0$ (as $P_{a,a}^0$), as well as $W^{a,b,c}:V^0\to V^a\otimes V^b\otimes V^c$ (as $(Id_a\otimes Y_a^{b,c})\circ\cup_a$) and $M_{a,b,c}:V^a\otimes V^b\otimes V^c\to V^0$ (as $\cap_c\circ(P_{a,b}^c\otimes Id_c)$) (see Figure \ref{fig:correspondence}). Composing these elementary morphisms one can associate morphisms to any diagram of a planar colored trivalent tangle.
\begin{prop}\label{prop:nonordered}
The following morphisms $V^a\otimes V^b\to V^c$ coincide:
$$\psfrag{=}{$=$}
\psfrag{d=}{$=$}
\psfrag{a}{\large$a$}
\psfrag{b}{\large$b$}
\psfrag{c}{\large$c$}
\psfrag{P}{\large$P_{a,b}^c=$}
\psfrag{rc}{\large$=(id_c\otimes \cap_b)\circ Y^{c,b}_a$}
\psfrag{lc}{\large$(\cap_a\otimes id_c)\circ Y^{a,c}_b$}
\includegraphics[width=6cm]{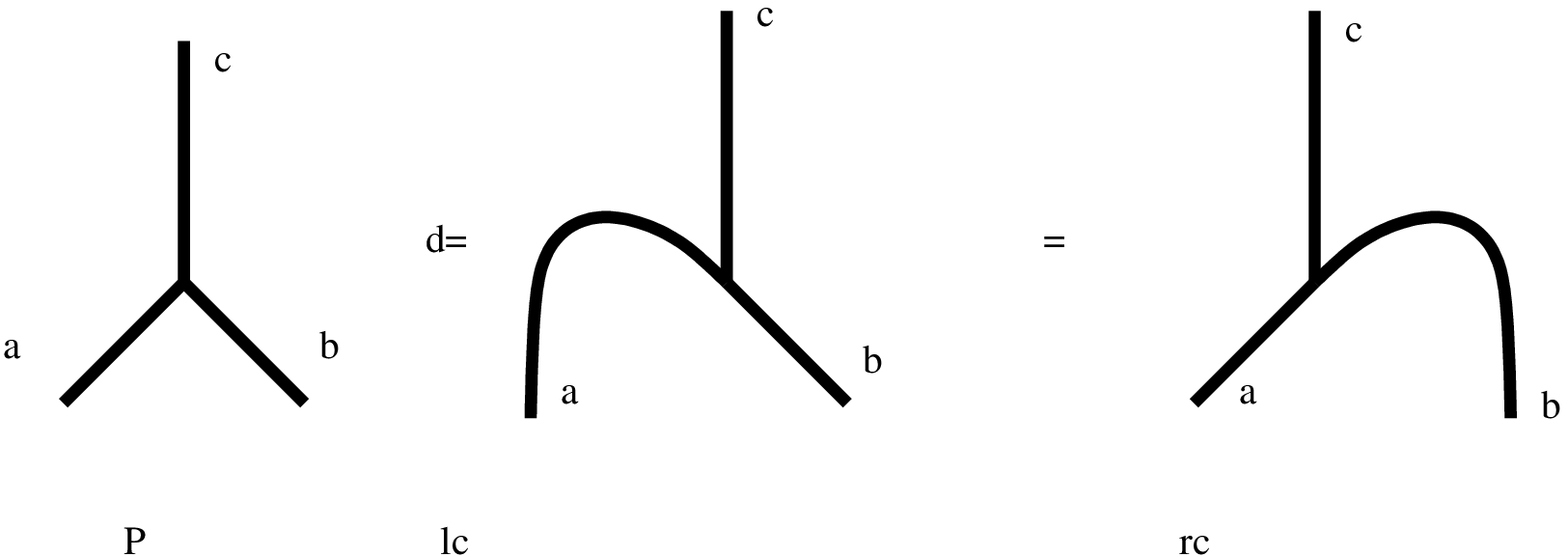}
$$\end{prop}
The proof is a direct consequence of isotopy invariance of morphisms induced by framed tangles colored by $\frac{1}{2}$, the definitions of $P^{c}_{a,b}$ and $Y^{a,b}_c$, the identities $JW_{2a}^2=JW_{2a}$ and the following: 
\begin{lemma}
The morphisms from $(V^{\frac{1}{2}})^{\otimes 4a}\to V^0$ represented by the following diagrams coincide:
$$\raisebox{-0.4cm}{\psfrag{2a}{\footnotesize{$2a$}}\includegraphics[width=1.5cm]{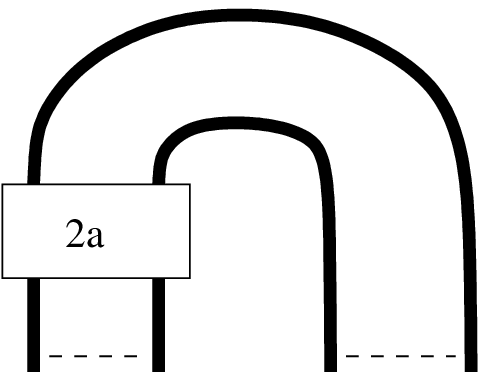}}=\raisebox{-0.4cm}{\psfrag{2a}{\footnotesize{$2a$}}\includegraphics[width=1.5cm]{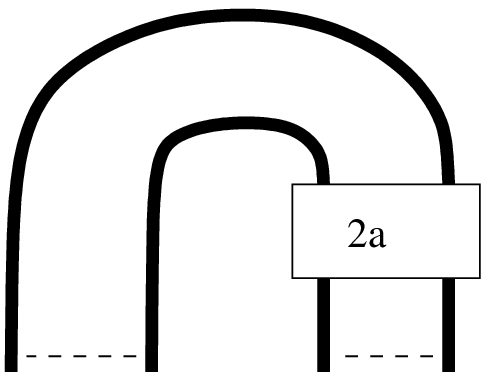}}$$
\end{lemma}
 \begin{prf}{1}{
Each minimal positive braid in the definition of $JW_{2a}$ of the l.h.s., can be slid through an isotopy over the max to a minimal positive braid in the definition of $JW_{2a}$ in the r.h.s. The morphisms two such braids induce are the same because the tangle they are represented by are isotopic, and their coefficients in the sum expressing $JW_{2a}$ are the same because they contain the same number of crossings. 
 }\end{prf}

One also defines operators ${}_{b}^{a}R$ (resp. ${}_{a}^bR_-$) $V^a\otimes V^b\to V^b\otimes V^a$ associated to a colored positive (resp. negative) crossing as:
$${}_{b}^{a}R\doteqdot (\mu_b\otimes\mu_a)\circ\raisebox{-0.7cm}{\psfrag{2a}{\footnotesize{$2a$}}\psfrag{2b}{\footnotesize{$2b$}}\psfrag{2c}{\footnotesize{$2c$}}\includegraphics[width=2.0cm]{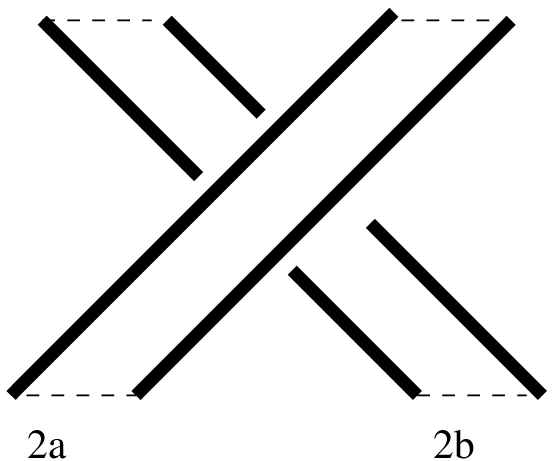}}\circ (\phi_a\otimes \phi_b)\ \ \ {\rm and}\ \ \ {}_{a}^{b}R_-\doteqdot (\mu_b\otimes\mu_a)\circ\raisebox{-0.7cm}{\psfrag{2b}{\footnotesize{$2a$}}\psfrag{2a}{\footnotesize{$2b$}}\psfrag{2c}{\footnotesize{$2c$}}\includegraphics[width=2.0cm]{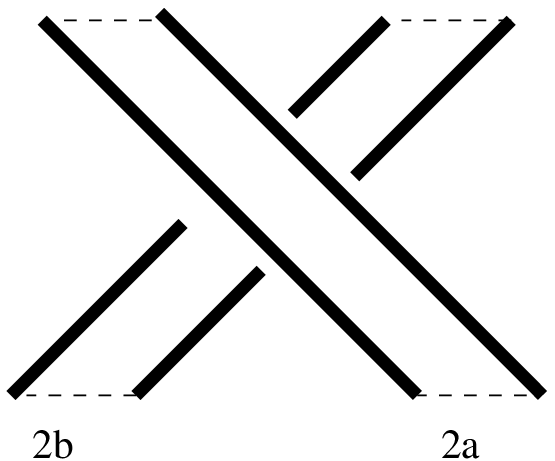}}\circ(\phi_a\otimes \phi_b)$$
where the diagrams represent $2a$ parallel strands passing over (resp. under) $2b$ parallel strands.

Now that we chose our elementary morphisms, let us conclude the proof of Proposition \ref{prop:kauffmanviauq}. 
Given a pair $(G,col)$ and a diagram $D$ for it, using the identities $\phi_a\circ\mu_a\circ JW_{2col(a)}=JW_{2col(a)}$ and  $JW^2_{2col(a)}=JW_{2col(a)}$ on all the edges of $D$, one sees that $op(G,col,D)=\sum_{i} c_i op(L_i,\frac{1}{2})$ where the $L_i$ and $c_i$ are the same framed links and coefficients as in the definition of $<G,col>$. But since we already proved that $op(L_i,\frac{1}{2})=<L_i,\frac{1}{2}>$, Proposition \ref{prop:kauffmanviauq} follows.

Later on, we will want to compute the invariants also using diagrams whose blackboard framing is not the one on $G$; therefore we define positive and negative ``half twist" endomorphisms $H_a:V^a\to V^a$ as:
$$H_a\doteqdot 
\mu_a\circ \raisebox{-0.6cm}{\psfrag{2a}{\footnotesize{$2a$}}\includegraphics[width=1.2cm]{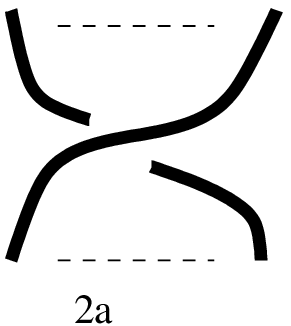}}\circ \phi_a \ \ {\rm and}\ \ H_a^{-1}\doteqdot 
\mu_a\circ \raisebox{-0.6cm}{\psfrag{2a}{\footnotesize{$2a$}}\includegraphics[width=1.2cm]{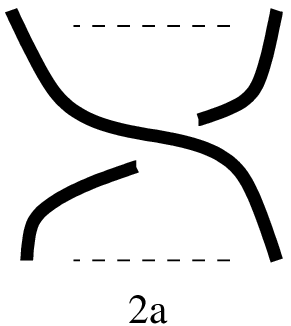}}\circ \phi_a$$
\subsection{Explicit formulas for the elementary operators}
In this subsection we provide explicit formulas for the operators defined in the preceding subsection when written in the bases $g^{a}_u$.
\subsubsection{Half-twists}
Let us start by the ``half twist" operator $H_a:V^a\to V^a$ which is $H_a\doteqdot(\sqrt{-1})^{2a} q^{(a^2+a)}Id_a$. We associate it to a vertical $a$-colored strand whose framing performs a positive half twist with respect to the blackboard framing. 
\subsubsection{Morphisms associated to $Y$-shaped vertices.}\label{sub:Y}
Let $Y_c^{a,b}\in Hom(V^c,V^a\otimes V^b)$ be defined as in Subsection \ref{sub:kauffmanUq}.
\begin{defi}
The \emph{quantum Clebsch-Gordan} coefficient $C^{a,b,c}_{u,v,t}$ is the coefficient in the sum:
$$Y_c^{a,b}(g^c_t)=\sum_{u+v=t}C^{a,b,c}_{u,v,t}\ g^a_u\otimes g^b_v$$
\end{defi}
It is clear that $C^{a,b,c}_{u,v,t}$ is  zero if $u+v\neq t$ because $\Delta(K)(g^a_u\otimes g^b_v)=q^{u+v}g^a_u\otimes g^b_v$ and the weight of a vector is preserved by a morphism. It holds:
\begin{prop}
\begin{equation}\label{eq:c}
C^{a,b,c}_{u,v,t}=\sqrt{-1}^{c-a-b}(-1)^{(b-v)-(c-t)}q^{\frac{(b-v)(b+v+1)-(a-u)(a+u+1)}{2}}\frac{[a+b-c]![b+c-a]![c+a-b]!}{[2c]!}\times$$
$$\times \sum_{z,w: z+w=c-t} (-1)^{z}q^{\frac{(z-w)(c+t+1)}{2}}\qbin{a+u+z}{a+c-b}\qbin{b+v+w}{b+c-a}\qbin{c-t}{z}
\end{equation}
where the sum is taken over all $z,w\in \mathbb{N}$ such that $z+w=c-t$ and all the arguments of the factorials are non-negative integers.
\end{prop}
\begin{prf}{1}{
In \cite{CS}, Lemma 3.6.10, using the basis $e^a_u=[a+u]!g^a_u$ for $V^a$, the following formula was provided for the Clebsch-Gordan coefficients (where we are rewriting the formula via $q$-binomials and correcting a missing factor of $\sqrt{-1}^{(t-c)}$):
\begin{equation}\label{eq:coldbase}
C^{a,b,c}_{u,v,t}=\sqrt{-1}^{c-a-b}(-1)^{(b-v)-(c-t)}q^{\frac{(b-v)(b+v+1)-(a-u)(a+u+1)}{2}}\frac{[c+t]![c-t]!}{[2c]!}\times$$
$$\times \sum_{z,w: z+w=c-t} (-1)^{z}q^{\frac{(z-w)(c+t+1)}{2}}\left[\begin{array}{c} a+b-c\\ a-u-z\end{array} \right]\left[\begin{array}{c} a+u+z\\ z\end{array}\right]\left[\begin{array}{c} b+v+w\\ w\end{array} \right]
\end{equation}
where the sum is taken over all $z,w$ such that $z+w=c-t$ and all the arguments of the quantum factorials are non-negative integers. 
To get our statement it is then sufficient to multiply by $\frac{[a+u]![b+v]!}{[c+t]!}$ (to operate the change of basis from $e^{a}_u$ to $g^{a}_u$), to single out of the factorials the terms $\frac{[a+b-c]![b+c-a]![c+a-b]!}{[2c]!}$ and to pair the factorials in the denominators of the summands so that their sums match $a+u+z$, $b+v+w$ and $c-t$ (recall that $u+v=t$).  
}\end{prf}

\subsubsection{Cup and Cap}
Let $\cup_a: \mc=V^a\to V^a\otimes V^a$ be defined as $Y_0^{a,a}$. An explicit computation using formula \ref{eq:c} gives in the base $g^a_u$:
\begin{equation}\label{eq:cup}
\cup_a(g^0_0)=\sum_{u=-a}^{a}[2a]!\sqrt{-1}^{2u}q^u g^a_u\otimes g^a_{-u}
\end{equation}
This, together with the invariance under isotopy which forces the identity $(\cap_a\otimes Id_a)\circ(Id_a\otimes \cup_a)=Id_a$ uniquely determines $\cap_a:V^a\otimes V^a\to V^0$ (defined as $P_{a,a}^0$) as:
\begin{equation}\label{eq:cap}
\cap_a(g^{a}_{u}\otimes g^{a}_v)=\delta_{u,-v}\frac{\sqrt{-1}^{2u}q^u}{[2a]!} g^0_0
\end{equation}
\subsubsection{Morphisms associated to $3$-valent vertices}\label{sub:P}
To compute the coefficients of the projectors $P_{a,b}^c:V^a\otimes V^b\to V^c$ out of $Y^{a,b}_c$ we use Proposition \ref{prop:nonordered}. So letting $P(g^a_{u}\otimes g^b_{v})=\sum_tP^{a,b,c}_{u,v,t} g^c_t$, it holds: 
\begin{equation}\label{eq:p}
P^{a,b,c}_{u,v,t}= C^{a,c,b}_{-u,t,v}\frac{\sqrt{-1}^{2u}q^{u}}{[2a]!}
\end{equation}

Similarly, ``non-smooth minima" operators $W^{a,b,c}\in Hom(V^0,V^a\otimes V^b\otimes V^c)$ are defined by ``pulling up the bottom leg in $Y_c^{a,b}$\ ", i.e. by setting $W^{a,b,c}\doteqdot (Y^{a,b}_c\otimes Id_c)\circ \cup_c$. So, letting: 
$$W^{a,b,c}(g^0_0)\doteqdot \sum_{u=-a}^{a}\sum_{v=-b}^{b}\sum_{t=-c}^{c} W^{a,b,c}_{u,v,t}\  g^a_u\otimes g^b_v\otimes g^c_t$$ the coefficients are:
\begin{equation}\label{eq:w}
W^{a,b,c}_{u,v,t}\doteqdot C^{a,b,c}_{u,v,-t} \sqrt{-1}^{-2t}q^{-t} [2c]!=
\sqrt{-1}^{-a-b-c}(-1)^{(b-v)-2t}q^{\frac{(b-v)(b+v+1)-(a-u)(a+u+1)}{2}}\times$$
$$\times [a+b-c]![b+c-a]![c+a-b]! \sum_{z,w: z+w=c+t} (-1)^{z}q^{-t+\frac{(z-w)(c-t+1)}{2}}\qbin{a+u+z}{a+c-b}\qbin{b+v+w}{b+c-a}\qbin{c+t}{z}
\end{equation}



Finally, $M^{a,b,c}\in Hom(V^a\otimes V^b\otimes V^c,V^0)$ defined by $M^{a,b,c}\doteqdot \cap_c\circ (P^{c}_{a,b}\otimes Id_c)$ has coefficients:
\begin{equation}\label{eq:m} 
M^{a,b,c}_{u,v,t}\doteqdot M^{a,b,c}(g^a_{u}\otimes g^b_{v}\otimes g^c_{t})=P^{a,b,c}_{u,v,-t}\frac{\sqrt{-1}^{-2t}q^{-t}}{[2c]!}
\end{equation}

\subsubsection{$R$-matrix.}
A positive crossing corresponds to the action of Drinfeld's universal $R$-matrix:
\begin{lemma}\label{lem:rmatrix}
The morphism ${}^{a}_{b}R:V^{a}\otimes V^{b}\to V^{b}\otimes V^{a}$ is given by the composition of Drinfeld's universal $R$-matrix with the flip of the coordinates and in the basis $g^a_u\otimes g^b_v$ it is:
\begin{equation}\label{eq:rmatrix}
R(g^{a}_{u}\otimes g^{b}_{v})=\sum_{n\geq 0} [n]!(q-q^{-1})^n\qbin{a-u}{n}\qbin{a+u+n}{n}q^{2uv-n(u-v)-\frac{n(n+1)}{2}}g^{b}_{v-n}\otimes g^{a}_{u+n}
\end{equation}
where the sum is taken over all the $n$ such that $\vert u+n\vert\leq a$ and $\vert v-n\vert\leq b$. 
We will denote ${}^{a}_{b}R_{u,v}^{h,k}$ the coefficient of $R(g^{a}_u\otimes g^{b}_v)$ with respect to $g^{b}_h\otimes g^{a}_k$.
\end{lemma}
\begin{prf}{1}{
The first statement is a direct consequence of the definition of ${}_{b}^{a}R$ and of the fact that ${}_{\frac{1}{2}}^\frac{1}{2}R$ coincides with the composition of the action of Drinfeld's $R$-matrix and the flip.  
To compute the explicit entries of ${}_{b}^{a}R$, we use the formulas provided in \cite{KM} (Corollary 2.32: recall that $\overline{t}=q^{-\frac{1}{2}}$) in the basis $f^j_m$ and the diagonal change of basis $f^j_m=\frac{[2j]!}{[j-m]!} g^j_m$:
$$R(g^{a}_{u}\otimes g^{b}_{v})=\sum_{n\geq 0} \frac{(q-q^{-1})^n}{[n]!} \frac{[a+u+n]![b-v+n]!}{[a+u]![b-v]!}\frac{[a-u]!}{[2a]!}\frac{[b-v]!}{[2b]!}\times$$ 
$$\times\frac{[2a]![2b]!}{[a-u-n]![b-v+n]!}q^{2uv-n(u-v)-\frac{n(n+1)}{2}}g^{b}_{v-n}\otimes g^{a}_{u+n} =$$
$$=\sum_{n\geq 0} [n]!(q-q^{-1})^n\qbin{a-u}{n}\qbin{a+u+n}{n}q^{2uv-n(u-v)-\frac{n(n+1)}{2}}g^{b}_{v-n}\otimes g^{a}_{u+n}$$}\end{prf}
The morphism associated to a negative crossing whose upper strand is colored by $a$ and whose lower strand is colored by $b$ is the inverse of ${}^{a}_{b}R$ and can be computed in terms of the one we just gave and two extrema: $(R_-)\doteqdot (Id_{a}\otimes Id_{b}\otimes \cap_{a})\circ(Id_{a}\otimes {}_{b}^{a}R\otimes Id_{a})\circ (\cup_{a}\otimes Id_{a}\otimes Id_{a})$.  An explicit formula is then computed out of formulas \ref{eq:rmatrix}), \ref{eq:cup}) and  \ref{eq:cap}):
\begin{equation}\label{eq:rmatrixinv}
{}_b^aR_-(g^{b}_{v}\otimes g^{a}_{u})=\sum_{n\geq 0}[n]! (q^{-1}-q)^{n}q^{-2vu+n(u-v)+\frac{n(n+1)}{2}}\qbin{a-u}{n}\qbin{a+u+n}{n} g^{a}_{u+n}\otimes g^{b}_{v-n}
\end{equation}
where the sum is taken over all the $n$ such that $\vert u+n\vert\leq a$ and $\vert v-n\vert\leq b$.
Remark that ${}_b^aR_-={}_b^aR\vert_{q\to q^{-1}}$ and that ${}_b^aR_-={}_b^aR^{-1}$  because ${}_b^aR_-\circ {}_b^aR=Id_a\otimes Id_b$. We will denote ${}^{a}_{b}(R^{-1})_{u,v}^{h,k}$ the coefficient of $R^{-1}(g^{a}_u\otimes g^{b}_v)$ with respect to $g^{b}_h\otimes g^{a}_k$.
The following well-known lemma relates $R$-matrices to $Y$-morphisms.
\begin{lemma}[Half-twist around a vertex]\label{lem:halftwist}
For every admissible $3$-uple $(a,b,c)$ it holds:
$$\raisebox{-.5cm}{\psfrag{a}{\footnotesize$a$}\psfrag{b}{\footnotesize$b$}\psfrag{c}{\footnotesize$c$}\includegraphics[width=1cm]{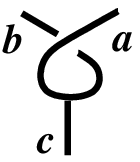}}={}_b^aR\circ Y^{a,b}_c=(H_b^{-1}\otimes H_a^{-1})\circ Y^{b,a}_c\circ(H^c)=\raisebox{-.4cm}{\psfrag{a}{\footnotesize$a$}\psfrag{b}{\footnotesize$b$}\psfrag{c}{\footnotesize$c$}\includegraphics[width=1.1cm]{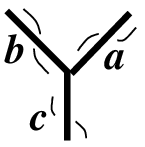}}$$
\end{lemma}
\begin{prf}{1}{
It is sufficient to prove the equality for the highest weight vector $g^c_c\in V^c$; to do this it is sufficient to check that ${}_b^aR\circ Y^{a,b}_c(g^c_c)$ and $(H_b^{-1}\otimes H_a^{-1})\circ Y^{b,a}_c\circ(H^c)(g^c_c)$ have the same coefficient with respect to the element $g^b_{b}\otimes g^a_{c-b}$ of the basis $g^b_j\otimes g^a_i,\ \vert i \vert\leq a,\ \vert j\vert \leq b$  of $V^b\otimes V^a$.
By formulas \ref{eq:rmatrix} and \ref{eq:c} this coefficient is for the left hand side:
$$\sum_{u+v=c} {}_b^aR^{b,c-b}_{u,v}\times C^{a,b,c}_{u,v,c}={}_b^aR^{b,c-b}_{c-b,b}\times C^{a,b,c}_{c-b,b,c}=$$
$$=q^{2(c-b)b}\times \sqrt{-1}^{c-a-b}(-1)^{0}q^{-\frac{1}{2}(a-(c-b))(a+(c-b)+1)}\frac{[2b]![a+c-b]!}{[2c]!}$$
The coefficient on the r.h.s. is $(H^b_b)^{-1} (H^a_a)^{-1} C^{b,a,c}_{b,c-b,c}(H^c_c)$ which equals:
$$\sqrt{-1}^{2c-2a-2b}q^{c^2+c-(a^2+a)-(b^2+b)}\sqrt{-1}^{c-b-a}(-1)^{a-(c-b)}q^{\frac{1}{2}(a-(c-b))(a+(c-b)+1)}\frac{[2b]![a+c-b]!}{[2c]!}=$$ $$=\sqrt{-1}^{c-b-a}q^{c^2+c-(a^2+a)-(b^2+b)}q^{\frac{1}{2}(a-(c-b))(a+(c-b)+1)}\frac{[2b]![a+c-b]!}{[2c]!}$$
A straightforward computation shows that the two coefficients are indeed equal.
}
\end{prf}

\subsection{The state-sum computing $\langle G,col\rangle$.}\label{sub:statesum}
Let $G$ be a closed KTG, $E$ be the set of its edges, $V$ the set of its vertices and $col:E\to \frac{\mathbb{N}}{2}$ an admissible coloring. 
Let also $D$ be a diagram of $G$ and for every $e\in E$, let $g_e\in \frac{\mathbb{N}}{2}$ be the difference between the framing of $e$ in $G$ and the blackboard framing on it (it is half integer because the two framings may differ of an odd number of half twists). 
Let $C,M,N$ be respectively the set of crossings, maxima and minima in $D$ (recall that we are fixing a height function on $\mr^2$ to decompose $D$ into elementary subgraphs). Then let $f_1,\ldots f_n$ be the connected components of $D\setminus (V\cup M\cup N\cup C)$. Remark that each $f_j$ is a substrand of an edge of $G$ therefore it inherits a color which we will denote $c_j$.
To express $\langle G,col\rangle$ as a state-sum, let us first define a state:
\begin{defi}[States]
A \emph{state} is a map $s:\{f_1,\ldots f_n\}\to \frac{\mz}{2}$ such that $\vert s(f_j)\vert \leq c_j,\  \forall j=1,\ldots n$ and $s(f_j)$ is integer iff $c_j$ is. Given a state $s$, we will call the value $s(f_j)$ the \emph{state of} $f_j$. (Equivalently, a state is a choice of one vector $g^{c_j}_{s(f_j)}$ of the base of $V^{c_j}$ for each component of $D\setminus (V\cup M\cup N\cup C)$.) 
\end{defi}
The \emph{weight} $w(s)$ of a state $s$ is the product of a factor $w_s(x)$ per each $x$ crossing, vertex, maximum and minimum of $D$. To define these factors, in the following table use the letters $a,b,c$ for the colors of the strands and $u,t,v,w$ for their states.
\vspace{.1cm}
\begin{center}\begin{tabular}{|c|c|}
\hline
$w_s(\raisebox{-.1cm}{\psfrag{a}{\footnotesize$a$}\psfrag{u}{\footnotesize$u$}\psfrag{v}{\footnotesize$v$}\includegraphics[width=1cm]{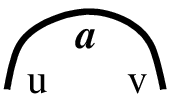}})=\delta_{u,-v}\frac{\sqrt{-1}^{2u}q^u}{[2a]!}$ & $w_s(\raisebox{-.2cm}{\psfrag{a}{\footnotesize$a$}\psfrag{u}{\footnotesize$u$}\psfrag{v}{\footnotesize$v$}\includegraphics[width=1cm]{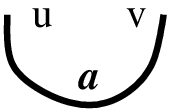}})=\delta_{v,-u} [2a]!\sqrt{-1}^{2u}q^u$\\
$w_s(\raisebox{-.3cm}{\psfrag{a}{\footnotesize$a$}\psfrag{b}{\footnotesize$b$}\psfrag{c}{\footnotesize$c$}\psfrag{u}{\footnotesize$u$}\psfrag{v}{\footnotesize$v$}\psfrag{t}{\footnotesize$t$}\includegraphics[width=.8cm]{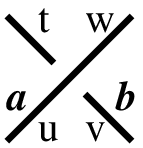}})={}_b^aR_{u,v}^{t,w}$ (see Formula \ref{eq:rmatrix}) & $w_s(\raisebox{-.3cm}{\psfrag{a}{\footnotesize$a$}\psfrag{b}{\footnotesize$b$}\psfrag{c}{\footnotesize$c$}\psfrag{u}{\footnotesize$u$}\psfrag{v}{\footnotesize$v$}\psfrag{t}{\footnotesize$t$}\includegraphics[width=.8cm]{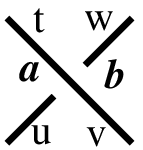}})={}_b^a(R^{-1})_{u,v}^{t,w}$ (see Formula \ref{eq:rmatrixinv})\\
 $w_s(\raisebox{-.5cm}{\psfrag{a}{\footnotesize$a$}\psfrag{b}{\footnotesize$b$}\psfrag{c}{\footnotesize$c$}\psfrag{u}{\footnotesize$u$}\psfrag{v}{\footnotesize$v$}\psfrag{t}{\footnotesize$t$}\includegraphics[width=0.8cm]{smallY.eps}})=C^{a,b,c}_{u,v,t}$ (see Formula \ref{eq:c}) & 
$w_s(\raisebox{-.4cm}{\psfrag{a}{\footnotesize$a$}\psfrag{b}{\footnotesize$b$}\psfrag{c}{\footnotesize$c$}\psfrag{u}{\footnotesize$u$}\psfrag{v}{\footnotesize$v$}\psfrag{t}{\footnotesize$t$}\includegraphics[width=0.8cm]{smallP.eps}})=P^{a,b,c}_{u,v,t}$ (see Formula \ref{eq:p})\\
 $w_s(\raisebox{-.2cm}{\psfrag{a}{\footnotesize$a$}\psfrag{b}{\footnotesize$b$}\psfrag{c}{\footnotesize$c$}\psfrag{u}{\footnotesize$u$}\psfrag{v}{\footnotesize$v$}\psfrag{t}{\footnotesize$t$}\includegraphics[width=1.2cm]{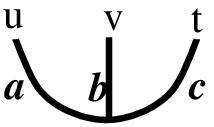}})=W^{a,b,c}_{u,v,t}$ (see Formula \ref{eq:w}) & $w_s(\raisebox{-.1cm}{\psfrag{a}{\footnotesize$a$}\psfrag{b}{\footnotesize$b$}\psfrag{c}{\footnotesize$c$}\psfrag{u}{\footnotesize$u$}\psfrag{v}{\footnotesize$v$}\psfrag{t}{\footnotesize$t$}\includegraphics[width=1.2cm]{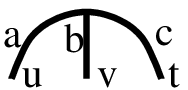}})=M^{a,b,c}_{u,v,t}$ (see Formula \ref{eq:m})\\
\hline
\end{tabular}
\end{center}
\vspace{.1cm}
Finally, to take into account the action of the half-twist operators $H_j$ on the edges of $G$, let
 $F(D,col)\doteqdot \prod_{e\in edges} \sqrt{-1}^{4g_e col(e)}q^{2g_e(col(e)^2+col(e))}$ be the \emph{framing factor} (note that it does not depend on any state but it does depend on the diagram $D$). The weight of a state $s$ of $\langle G,col\rangle$ is then defined as follows:
\begin{equation}\label{eq:stateweight}
w(s)=F(D,col)\prod_{M\in maxima}w_s(M)\prod_{m\in minima} w_s(m)\prod_{c\in crossings}w_s\large(c)\prod_{v\in vertices} w_s(v)
\end{equation}
The value of $\langle  G,col\rangle $ is then given by :
\begin{equation}\label{eq:statesum}
\langle  G,col\rangle =\sum_{s\in states} w(s)
\end{equation}
since the state-sum represents nothing else than the composition of the elementary morphisms associated to $G$ as a morphism $op(G,col,D):V^0\to V^0$.
\begin{rem}
The above state-sum shows that one can extend this way the definition of Kauffman brackets to colored KTG's whose framing is a non-orientable surface: diagrams of such KTG's will always contain some half-twists which will contribute through a constant multiplicative factor (included in $F(D,col)$). 
\end{rem}
\begin{example}\label{ex:unlinks}
If $L=L_1,\ldots L_n$ is an unlink with (possibly half-integral) framings $g_1,\ldots g_n$ and colored by colors $c_1,\ldots c_n$ then:
$$\langle  L,col\rangle =\prod_{i=1}^{n} (-1)^{2c_i}[2c_i+1](\sqrt{-1})^{4g_ic_i}q^{2g_i(c_i^2+c_i)}$$
It is sufficient to prove it for the case of an unknot colored by $c_j$ and with framing $g_j$.
In that case the value is the trace of the operator $\cap_{c_j}\circ (Id\otimes H_{c_j}^{2g_j})\circ \cup_{c_j}$ which equals: 
$$tr(\cap_{c_j}\circ (Id\otimes H_{c_j}^{2g_j})\circ \cup_{c_j})=(\sqrt{-1})^{4g_ic_i}q^{2g_i(c_i^2+c_i)}\sum_{u=-c_j}^{c_j} \sqrt{-1}^{4u}q^{2u}=$$
$$(\sqrt{-1})^{4g_ic_i}q^{2g_i(c_i^2+c_i)}(-1)^{2c_j}\sum_{u=-c_j}^{c_j} q^{2u}=(-1)^{2c_j}[2c_j+1](\sqrt{-1})^{4g_ic_i}q^{2g_i(c_i^2+c_i)}$$
\end{example}
\begin{figure}
\psfrag{a}{\large$a$}
\psfrag{b}{\large$b$}
\psfrag{c}{\large$c$}
\psfrag{=}{\large$=$}
\psfrag{dc}{\large$\doteqdot$}
\psfrag{M}{\large$M^{a,b,c}\circ W^{a,b,c}$}
\psfrag{sum}{\large$\sum W^{a,b,c}_{m_1,m_2,m}M^{a,b,c}_{m_1,m_2,m}$}

\includegraphics[width=5cm]{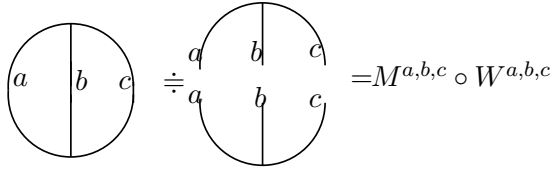}
\caption{The operator associated to a colored theta graph.}\label{fig:theta2}
\end{figure}
\begin{example}[The theta graph]
Consider a theta graph as in Figure \ref{fig:theta2}. Equip it with the blackboard framing and color the edges by $V^a,V^b$ and $V^c$. Its invariant is then: 
\begin{equation}\label{eq:theta2}
\theta(a,b,c)\doteqdot \sum_{m_1=-a}^a\sum_{m_2=-b}^b \sum_{m=-c}^c W^{a,b,c}_{m_1,m_2,m}M^{a,b,c}_{m_1,m_2,m}
\end{equation}
By Proposition \ref{prop:kauffmanviauq} the above formula equals the standard skein theoretical value:
\begin{equation}\label{eq:theta}
\theta(a,b,c)=(-1)^{a+b+c}\frac{[a+b+c+1]![a+b-c]![a+c-b]![b+c-a]!}{[2a]![2b]![2c]!}
\end{equation}
\end{example}
\begin{rem}
The above example shows that in general $\langle  G,col\rangle $ is not a Laurent polynomial: consider for instance the case $a=b=c=1$ in Formula \ref{eq:theta}.
\end{rem}

\section{Integrality}\label{sec:main}
Let $G$ be a closed KTG, $E$ be the set of its edges, $V$ the set of its vertices and $col:E\to \frac{\mathbb{N}}{2}$ an admissible coloring. We define the \emph{Euler characteristic} $\chi(e)$ of an edge $e\in E$ as $1$ if $e$ touches a vertex and $0$ otherwise (some edges of $G$ may be knots).
Let also $D$ be a diagram of $G$ and for every $e\in E$, let $g_e\in \frac{\mathbb{N}}{2}$ be the difference between the framing of $e$ in $G$ and the blackboard framing on it (it is half integer because the two framings may differ of an odd number of half twists). 
We define a renormalization for $\langle G,col\rangle$ as follows:
$$
\langle  \langle  G,col\rangle \rangle \doteqdot \langle  G,col\rangle  \frac{\prod_{e\in E, \chi(e)=1}[2col(e)]!}{\prod_{v\in V}[a_v+b_v-c_v]![b_v+c_v-a_v]![c_v+a_v-b_v]!}
$$
where by $a_v,b_v,c_v$ we denote the colors of the three edges surrounding $v$.
\begin{rem}
The renormalization factor depends only on the abstract combinatorial structure of $(G,col)$, therefore $<<G,col>>$ is an invariant of colored KTG's.
\end{rem}
\begin{teo}[Integrality of the renormalized Kauffman brackets]\label{teo:main}
There exist $m,n\in \mz$ such that:
$$
\langle  \langle  G,col\rangle \rangle  \in (\sqrt{-1})^mq^{\frac{n}{4}}\mathbb{Z}[q,q^{-1}]
$$
Moreover, $<<G,col>>(\sqrt{-1})^{-m}q^{-\frac{n}{4}}$ is divisible in $\mathbb{Z}[q,q^{-1}]$ by $[2col(e)+1]$ for each edge $e$ of $G$. If the framing of $G$ is orientable, then $m=0$ and $n$ is even.
\end{teo}
\begin{prf}{1}{
Up to isotopy, we can suppose that the diagram $D$ of $G$ is the closure of a $(1,1)$-tangle $G'$ whose boundary strands are included in $e$ and also (by small isotopies around vertices and crossings) that $D$ contains only maxima, minima, positive crossings and vertices with $3$ top legs. The factor $F(D,col)=\sqrt{-1}^{\sum_{e\in E}-4g_e col(e)}\prod_{e\in E}q^{-2g_e(col(e)^2+col(e))}$ in the state-sum \ref{eq:stateweight}) changes the value of $<G,col>$ only by a factor of the form $(\sqrt{-1})^{k}q^{\frac{h}{4}},\ k,h\in \mz$, therefore, up to dividing by $F(D,col)$ we may suppose that the framing of $G$ is the blackboard framing. By Schur's lemma the morphism represented by $G'$ is $\lambda Id_{V^n}$. We claim that there exists an integer $s$ such that $\mu\doteqdot \lambda \frac{\prod_{e\in E, \chi(e)=1}[2col(e)]!}{\prod_{v\in V}[a_v+b_v-c_v]![b_v+c_v-a_v]![c_v+a_v-b_v]!}$ belongs to $q^{\frac{s}{2}}\mathbb{Z}[q,q^{-1}]$; this will conclude because $\langle\langle G,col\rangle\rangle=(-1)^{2col(e)}[2col(e)+1]\mu$.

To prove our claim let us define ``renormalized operators" associated to each maximum, minimum, crossing and vertex of $D$ equipped with a state as follows:
$$NH_j\doteqdot H_j\ \ \ \ \ \ \ (N\cup_a)_{u,v}\doteqdot \delta_{u,-v}\frac{1}{[2a]!}\raisebox{-.2cm}{\includegraphics[width=1cm]{smallmin.eps}}\ \ \ \ \ \ (N\cap_a)_{u,v}\doteqdot\delta_{u,-v}[2a]!\raisebox{-.2cm}{\includegraphics[width=1cm]{smallmax.eps}}$$
$$({}^{a}_{b}NR)^{t,w}_{u,v}:={}^{a}_{b}R^{t,w}_{u,v}\ \ \ \ \ \ (NW)^{a,b,c}_{u,v,t}\doteqdot \frac{W^{a,b,c}_{u,v,t}}{[a+b-c]![b+c-a]![c+a-b]!} $$

Since the morphism represented by $G'$ is diagonal, the only non-zero weight states are those where the states of the top and bottom strand of $G'$ are equal. Therefore, if in formula \ref{eq:statesum}) one fixes the same state $u$ on the top and bottom strand of $G'$ and replaces each weight by its ``normalized version" defined above, the result will be:
$$
\langle\langle G',col\rangle\rangle =\langle  G',col\rangle \times \frac{\prod_{e\in E}([2col(e)]!)^{(cap(e)-cup(e))}}{\prod_{v\in V}[a_v+b_v-c_v]![b_v+c_v-a_v]![c_v+a_v-b_v]!}
$$ where for each edge, $cup(e)$ (resp. $cap(e)$) are the number of minima (resp. maxima) on $e$. The above formula coincides with normalization factor as in the claim since by our hypothesis on $D$ all the vertices have $3$ top legs and so for each edge $e$ different from the top strand it holds $\chi(e)=cap(e)-cup(e)$, and that for the top strand $cap(e)-cup(e)=0$ (the top and bottom strands are part of the same edge in $G$ therefore only one of them should be counted in the renormalization). 

Remark now that that each ``renormalized operator" has coefficients in $\mz[\sqrt{-1}][q^{\pm\frac{1}{2}}]$. This is straightforward because of Lemma \ref{lem:integrality} and Formulas \ref{eq:rmatrix} and \ref{eq:w}. 
Let us first show that actually all the coefficients are non-imaginary. Let us remark that in the state-sum for $<<G,col>>$, since each edge with $\chi(e)=1$ has two endpoints and each $NW^{a,b,c}$ belongs in particular to $\sqrt{-1}^{-a-b-c}\mz[q^{\pm\frac{1}{2}}]$, the product of the factors $\sqrt{-1}^{col(e)}$ coming from these vertices is $\sqrt{-1}^{-2col(e)}$. Similarly the product of the factors $\sqrt{-1}^{2u_i}$ coming from the cups and caps on $e$ is $\pm\sqrt{-1}^{2col(e)(cap(e)-cup(e))}=\pm\sqrt{-1}^{2col(e)}$ (because each $u_i$ is half integer iff $col(e_i)$ is) and this cancels with the previous imaginary phase.

So now we are left to show that $<<G,col>>$ contains only all odd or all even powers of $q^{\frac{1}{2}}$. First of all remark that this is the case for the coefficients $NW^{a,b,c}_{u,v,t},({}_b^aNR)_{u,v}^{v-n,u+n},(N\cap_a)_{u,v},(N\cup_a)_{u,v}$; more specifically, an inspection in Formulas  \ref{eq:w} and \ref{eq:rmatrix}, reveals that the parities ($\in \frac{\mz}{2\mz}$) of the powers of $q^{\frac{1}{2}}$ in these coefficients are (beware: states may be half-integers but the values below are integers, then considered mod $2$): 
\begin{itemize}
\item in $NW^{a,b,c}_{u,v,t}$: $a(a+1)-u(u+1)+b(b+1)-v(v+1)+c(c+1)-t(t+1)$;
\item in $({}_b^aNR)_{u,v}^{v-n,u+n}$: $i(a)n+i(b)n+i(a)i(b)$ where $i(x)\doteqdot1$ if $2x$ is odd and $0$ otherwise;
\item in $(N\cup_a)_{u,v}$ and $(N\cap_a)_{u,v}$: $2a$ (not depending on any state). 
\end{itemize}
Hence for each state $s$ in the state-sum expressing $<<G,col>>$ the weight $w(s)$ contains only even or odd powers of $q$; we will therefore call $s$ \emph{even} or \emph{odd} accordingly.
Our goal is to show that all the states have the same parity: for instance remark that if $col$ has integers values, all the states are even. 
Now, for each state $s$ we will compute its parity by ``redistributing" on the edges the parities of the coefficients of the elementary operators and then summing them up over all the edges. For each edge $e$ of $G$ with $\chi(e)=1$, orient $e$ arbitrarily and let $u$ and $v$ the states of the substrands of $e$ respectively at the beginning and at the end of $e$; let the contributions of the endpoints  to the \emph{parity} of $e$ be $2(col(e)^2+col(e))-u^2-u-v^2-v$. Similarly each $\cup$ or $\cap$ in $e$ contributes by $2col(e)$ and since $e$ contains an odd number of such operators they contribute globally by $2col(e)$. Finally to take into account the crossings, follow $e$ and remark that each time $e$ crosses another edge, say $e'$, the state on the substrand of $e$ jumps from $x$ to $x\pm n$ ($n\in \mz$) and the parity of the powers of $q^{\frac{1}{2}}$ in the $R$-matrix corresponding to the crossing is $ni(col(e))+ni(col(e'))+ i(col(e))i(col(e'))$; so we define the contribution of the crossing to the parity of $e$ as $ni(col(e))$, dropping for the moment the term $i(col(e))i(col(e'))$ which does not depend on the state. Using the fact that on each $\cup_{col(e)}$ and $\cap_{col(e)}$ the state of the substrand of $e$ changes sign (see Formulas \ref{eq:cup},\ref{eq:cap}) and following $e$ from its beginning to its end, one can check that the global contribution to the parity of $e$ coming from the crossings is $-(u+v)i(col(e))$. Therefore summing up the parity of $e$ is $2(col(e)^2+col(e))-u^2-u-v^2-v+2col(e)-(u+v)i(col(e))$, which modulo $2$ is $2col(e)^2-u^2-v^2$ if $i(col(e))=1$ or $0$ if $i(col(e))=0$. In both cases it is constant mod $2$ when $u$ and $v$ range in $\{-col(e),-col(e)+1,\ldots col(e)\}\subset \frac{\mathbb{Z}}{2}$. Similarly, for edges with $\chi(e)=0$ the parity is easily seen to be constant. Therefore the parity of the states is constant because it is the sum of constant contributions coming from the edges and the constant term $C=\sum_{crossings} i(a)i(b)$ (where $a$ and $b$ are the colors of the strands forming the crossing).

To prove the last statement remark that by Lemma \ref{lem:orientableframing} one can choose $D$ satisfying the additional requirement that the framing of $G$ coincides with the blackboard framing, and so in this case one can suppose $F(D,col)=1$. In the above proof, this factor was the only one responsible for possible terms of the form $\sqrt{-1}q^{\frac{1}{4}}$.
}\end{prf}
\begin{cor}
With the notation of Theorem \ref{teo:main}, if $a_1,\ldots a_k$ are colors of edges of $G$ such that $2a_i+1,i=1,\ldots k$ are pairwise co-prime, $\langle\langle G,col\rangle\rangle\sqrt{-1}^{-m}q^{\frac{-n}{4}}$ is divisible by $\prod_j [2a_j+1]$ in $\mathbb{Z}[q^{\pm\frac{1}{2}}]$.
\end{cor}
\begin{lemma}\label{lem:integrality}
Let $a_1,\ldots a_s$ be integers and let the q-multinomial be defined as 
$$ \left[ \begin{array}{c}
a_1+\cdots +a_s \\
a_1\ \ a_2\ \  \ldots \ \ a_{s-1}\ \  a_s
\end{array}\right]\doteqdot \frac{[a_1+\cdots +a_s]!}{[a_1]!\cdots [a_s]!}
$$ 
Then the $q$-multinomial is a Laurent polynomial with positive, integer coefficients.
\end{lemma}
\begin{prf}{1}{
If $s=2$ the statement is a direct consequence of the fact that, if $yx=q^2xy$ are two skew-commuting variables, then: 
$$(x+y)^n=\sum_{j=0}^n q^{\frac{n(n+1)}{2}-\frac{j(j+1)}{2}-\frac{(n-j)(n-j+1)}{2}}\left[ \begin{array}{c}
n\\
j
\end{array}\right] x^jy^{n-j}
$$
which is easily proved by induction. The general case follows by induction on $s$ by remarking that:
$$ \left[ \begin{array}{c}
a_1+\cdots +a_s \\
a_1\ \ a_2\ \  \ldots \ \ a_{s-1}\ \  a_s
\end{array}\right]= \left[ \begin{array}{c}
a_1+\cdots +a_s \\
(a_1+a_2)\ \  \ldots \ \ a_{s-1}\ \  a_s
\end{array}\right]\left[ \begin{array}{c}
a_1+a_2 \\
a_1\ \ a_2
\end{array}\right]
$$ 
}\end{prf}
\subsection{Examples and properties}\label{sec:properties}
The following examples can be proved by re-normalizing the formulas provided in \cite{MV} for the standard skein invariants of tetrahedra and $\theta$-graphs.
\begin{example}[Unknot]
\begin{equation}\label{eq:nunknot}
\langle\langle\smallunknot{a}{.9}{.3}\rangle\rangle=(-1)^{2a}[2a+1]\end{equation}
\end{example}
\begin{example}[Theta graph]
\begin{equation}\label{eq:ntheta}
\langle  \langle \smalltheta{a}{b}{c}{.9}{.3}\rangle \rangle =(-1)^{a+b+c}[a+b+c+1]\left[\begin{array}{c}
a+b+c\\
a+b-c, b+c-a, c+a-b \end{array}\right]
\end{equation}
\end{example}
\begin{example}[The tetrahedron or symmetric $6j$-symbol]\label{ex:tet}
\begin{equation}\label{eq:ntet}
\langle  \langle\smalltet{a}{b}{c}{d}{e}{f}{.9}{.3}\rangle \rangle =\sum_{z= Max T_i}^{z=Min Q_j} (-1)^{z}[z+1]\left[\begin{array}{c}
z\\
z-T_1, z-T_2, z-T_3, z-T_4, Q_1-z, Q_2-z, Q_3-z \end{array}\right]
\end{equation}
where $T_1=a+b+c,\ T_2=a+e+f,\ T_3=d+b+f,\ T_4=d+e+c,\  Q_1=a+b+d+e,\ Q_2=a+c+d+f, \ Q_3=b+c+e+f$. 
\end{example}
\begin{example}[The crossed tetrahedron]\label{ex:crossedtet}
\begin{equation}\label{eq:ntetcross}
\smalltetcross{a}{b}{c}{d}{e}{f}{.9}{.3}=\raisebox{-.5cm}{\psfrag{a}{\footnotesize{$a$}}\psfrag{b}{\footnotesize{$b$}}\psfrag{c}{\footnotesize{$c$}}\psfrag{d}{\footnotesize{$d$}}\psfrag{e}{\footnotesize{$e$}}\psfrag{f}{\footnotesize{$f$}}\includegraphics[width=1.2cm]{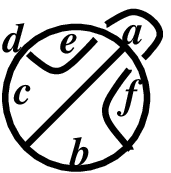}}=\raisebox{-.5cm}{\psfrag{a}{\footnotesize{$a$}}\psfrag{b}{\footnotesize{$b$}}\psfrag{c}{\footnotesize{$c$}}\psfrag{d}{\footnotesize{$d$}}\psfrag{e}{\footnotesize{$e$}}\psfrag{f}{\footnotesize{$f$}}\includegraphics[width=1.2cm]{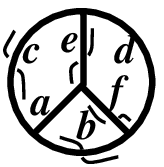}}=\sqrt{-1}^{2(f+c-e-b)}q^{f^2+f+c^2+c-b^2-b-g^2-e}\ \ \smalltet{e}{f}{a}{b}{c}{d}{.9}{.3}
\end{equation}
\end{example}

From now on we will often drop the $\langle\langle \cdot \rangle\rangle$ and denote the values provided by the above formulas respectively by: $$\smallunknot{a}{.9}{.3},\smalltheta{a}{b}{c}{.9}{.3},\smalltet{a}{b}{c}{d}{e}{f}{.9}{.3} \ \ \  {\rm and}\ \ \smalltetcross{a}{b}{c}{d}{e}{f}{.9}{.3} $$
All the invariants we will be dealing with will be the normalized version unless explicitly stated the contrary. 
The following lemma can be easily proved by starting from the analogous statements for the standard skein theoretical normalization of the invariants:
\begin{lemma}[Properties of the renormalized invariant]\label{lem:properties}
The following are some of the properties of $\langle  \langle  G,col\rangle \rangle $:
\begin{enumerate}
\item (Erasing $0$-colored strand) If $G'$ is obtained from $(G,col)$ by deleting a $0$-colored edge, then $\langle  \langle  G',col'\rangle \rangle =\langle  \langle  G,col\rangle \rangle $. 
\item (Connected sum) If $G=G_1\# G_2$ along an edge colored by $a$, then 
$$\langle  \langle  \raisebox{-.5cm}{\psfrag{a}{$a$}\psfrag{G1}{$G_1$}\psfrag{G2}{$G_2$}\includegraphics[width=3cm]{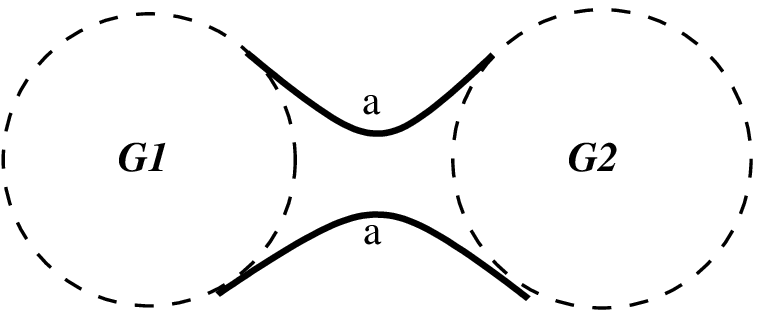}}\rangle \rangle =\frac{1}{(-1)^{2a}[2a+1]}\ \langle  \langle   \raisebox{-.5cm}{\psfrag{a}{$a$}\psfrag{G1}{$G_1$}\psfrag{G2}{$G_2$}\includegraphics[width=3cm]{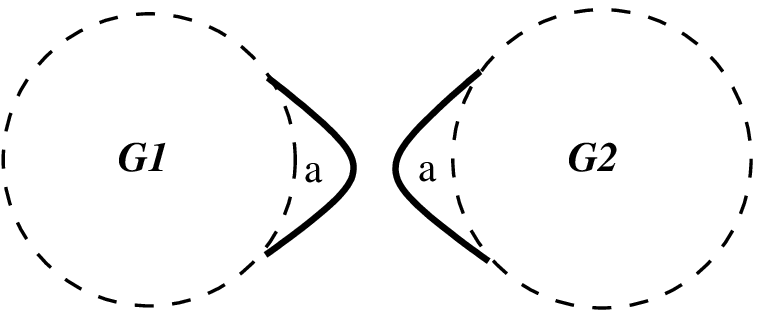}}\rangle \rangle $$
\item (Whitehead move) If $G$ and $G'$ differ by a Whitehead move then: $$\langle  \langle  \raisebox{-.6cm}{\psfrag{a}{$a$} \psfrag{b}{$b$} \psfrag{c}{$c$} \psfrag{d}{$d$} \psfrag{j}{$j$}\includegraphics[width=2cm]{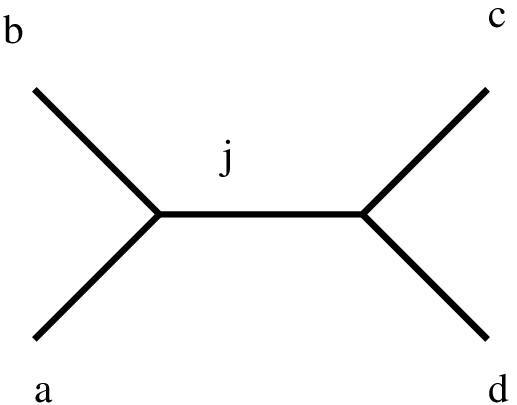}}\rangle \rangle =\Large{\sum}_{i} \frac{\smallunknot{i}{.9}{.3}\smalltet{a}{b}{i}{c}{d}{j}{.9}{.3}}{\smalltheta{a}{i}{d}{.9}{.3}\smalltheta{b}{i}{c}{.9}{.3}}\ 
\langle  \langle \raisebox{-1.0cm}{\psfrag{a}{$a$} \psfrag{b}{$b$} \psfrag{c}{$c$} \psfrag{d}{$d$} \psfrag{i}{$i$}\includegraphics[width=1cm]{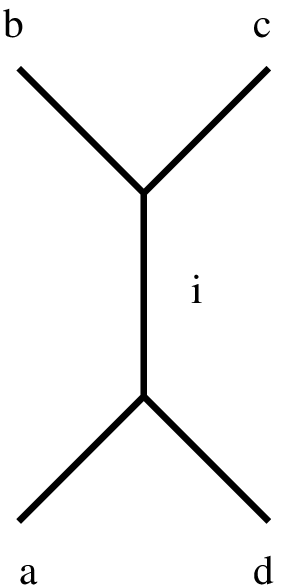}}\rangle \rangle $$ where $i$ ranges over the all the admissible values. 
\item (Fusion rule) In particular, applying the preceding formula to the case $j=0$ one has:
$$\langle  \langle  \raisebox{-.6cm}{\psfrag{a}{$a$} \psfrag{b}{$b$} \psfrag{c}{$c$} \psfrag{d}{$d$} \psfrag{j}{$j$}\includegraphics[width=1.cm]{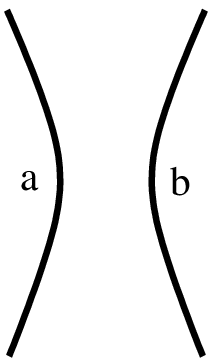}}\rangle \rangle =\Large{\sum}_{i}\frac{\smallunknot{i}{.9}{.3}}{\smalltheta{a}{i}{c}{.9}{.3}}\ \langle  \langle  \raisebox{-1.0cm}{\psfrag{a}{$a$} \psfrag{b}{$b$} \psfrag{c}{$b$} \psfrag{d}{$d$} \psfrag{i}{$i$}\includegraphics[width=1cm]{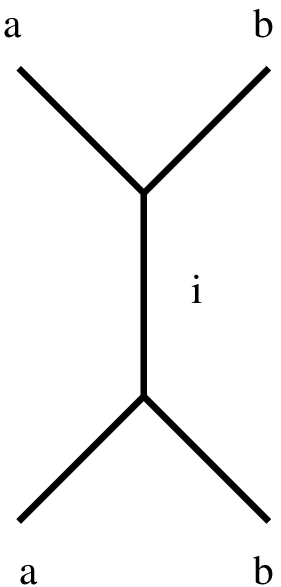}}\rangle \rangle $$
\end{enumerate}
\end{lemma}
\section{Shadow state-sums and integrality}\label{sec:shadows}
In this section we will first provide a so-called shadow-state sum formula for the invariants $\langle\langle G,col\rangle\rangle$ and give some examples. Then we will show through explicit examples that the summands of the shadow state-sums are not Laurent polynomials (even though, of course, the global sum of the state-sums are) in Subsection \ref{sub:commentsintegrality} will provide short proofs of known identities on $6j$-symbols. 
\subsection{Shadow state sums}
Let $(G,col)$ be a fixed colored graph, $D\subset \mr^2$ be a diagram of it and $V,E$ be the sets of vertices and edges of $G$, and $C,F$ the sets of crossings and edges of $D$ (each edge of $D$ is a sub-arc of one of $G$ therefore it inherits the coloring from $col$). Let the \emph{regions} $r_0,\ldots r_m$ of $D$ be the connected components of $\mr^2\setminus D$ with $r_0$ the unbounded one; we will denote by $R$ the set of regions and we will say that a region ``contains" an edge of $D$ or a crossing if its closure does. 

\begin{figure}
\includegraphics[width=10cm]{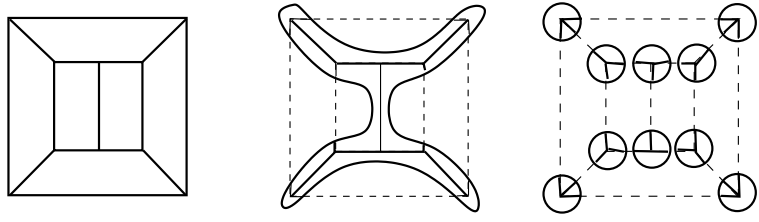}
\caption{On the left a planar graph $G$ to which we apply the construction of the proof of Theorem \ref{teo:shadowstate}. In the middle the graph $A$ constructed by shrinking on a maximal subtree of $G$ an unknot bounding a disk containing $G$ (the dotted parts are left just for reference). On the right the final union of planar tetrahedra (in this case there are no crossings).}\label{fig:planargraph}
\end{figure}

\begin{defi}[Shadow-state]
A shadow-state $s$ is a map $s:R\cup F \to \frac{\mathbb{N}}{2}$ such that $s(f)$ equals the color of the edge of $G$ containing $f$, $s(r_0)=0$ and  whenever two regions $r_i$ and $r_j$ contain an edge $f$ of $D$ then $s(r_i),s(r_j),s(f)$ form an admissible three-uple.
\end{defi}
Given a shadow-state $s$, we can define its weight as a product of factors coming from the local building blocks of $D$ i.e. the regions, the edges of $D$, the vertices of $G$ and the crossings. To define these factors explicitly, in the following we will denote by $a,b,c$ the colors of the edges of $G$ (or of $D$) and by $u,v,t,w$ the shadow-states of the regions and will use the examples given in Section \ref{sec:properties}.
\begin{enumerate}
\item If $r$ is a region whose shadow-state is $u$ a	nd $\chi(r)$ is its Euler characteristic, $$w_s(r)\doteqdot \smallunknot{u}{.9}{.4}^{\chi(r)}$$
\item If $f$ is an edge of $D$ colored by $a$ and $u,v$ are the shadow-states of the regions containing it then, letting $\chi(f)$ to be $0$ if $f$ is a closed component and $1$ otherwise: $$w_s(f)\doteqdot \smalltheta{u}{a}{v}{.9}{.4}^{\chi(f)}$$
\item If $v$ is a vertex of $G$ colored by $a,b,c$ and $u,v,t$ are the shadow-states of the regions containing it then $$w_s(v)\doteqdot \smalltet{a}{b}{c}{u}{v}{t}{.9}{.4}$$
\item If $c$ is a crossing between two edges of $G$ colored by $a,b$ and $u,v,t,w$ are the shadow-states of the regions surrounding $c$ then $$w_s(c)\doteqdot \smalltetcross{a}{u}{v}{b}{t}{w}{.9}{.4}$$
\end{enumerate}  
From now on, to avoid a cumbersome notation, given a shadow-state $s$ we will not explicit the colors of the edges of each graph providing the weight of the local building blocks of $D$ as they are completely specified by the states of the regions and the colors of the edges of $G$ surrounding the block. 
Then we may define the weight of the shadow-state $s$ as:
\begin{equation}\label{eq:shadowstate}
w(s)=\prod_{r\in R}  \smallunknot{}{.9}{.4}^{\chi(r)}\prod_{f\in F} \smalltheta{}{}{}{.9}{.4}^{-\chi(f)}\nns\nns\nns\prod_{v\in V}\nns \smalltet{}{}{}{}{}{}{.9}{.4}\prod_{c\in C}\smalltetcross{}{}{}{}{}{}{.9}{.4}
\end{equation}
Then, since the set of shadow-states of $D$ is easily seen to be finite, we may define the shadow-state sum of $(G,col)$ as
$$shs(G,col)\doteqdot\sum_{s\in shadow\ states} w(s)$$
As the following theorem says, the shadow state-sums provide a different approach to the computation of $\langle\langle G,col\rangle\rangle$:
\begin{teo}\label{teo:shadowstate}
It holds: $\langle\langle G,col\rangle\rangle=F(D,col)shs(G,col)$, where (as in the preceding sections)  $F(D,col)\doteqdot \prod_{e\in edges} \sqrt{-1}^{4g_e col(e)}q^{2g_e(col(e)^2+col(e))}$.
\end{teo}
The original definition of the shadow state sums and proof of the above result (but for the standard normalization of the invariants) is due to Kirillov and Reshetikhin (\cite{KR}) and was later generalized to general shadows by Turaev (\cite{Tu}). We used this formulation to extend the definition of colored Jones polynomials to the case of graphs and links in connected sums of copies of $S^2\times S^1$ (\cite{Co2}) and to prove a version of the generalized volume conjecture for an infinite family of hyperbolic links called fundamental hyperbolic links (\cite{Co}).  These links were already studied in \cite{CT} for their remarkable topological and geometrical properties.
\begin{prf}{1}{
Multiplying by $F(D,col)^{-1}$ we can reduce to the case when the framing of $G$ is the blackboard framing. Let $D$ be a diagram of $G$; we can add to $G$ some $0$-colored edges cutting the regions of $D$ (except $r_0$) into discs (this changes $G$ and $D$ but not the value of the resulting invariant by Lemma \ref{lem:properties}); for each region we will need $\chi(r)-1$ such arcs. Fix also a maximal connected sub-tree $T$ in $D$ and let $o\subset \mr^2$ be a $0$-colored unknot bounding a disc containing $D$; it is clear that $\langle\langle G,col\rangle\rangle=\langle\langle (G,col)\cup (o,0)\rangle\rangle$. Let also $A$ be the trivalent graph defined as follows: $A\doteqdot (N(T)\cap D)\cup \partial N(T)$ (where $N(T)$ is the regular neighborhood of $T$ in $\mr^2$).
The idea of the proof is to apply a sequence of fusion rules and inverse connected sums in order to express $\langle\langle G,col\rangle\rangle$ as a $\sum_{col_i} c(col_i) \langle\langle A,col_i\rangle\rangle$ for some colorings $col_i$ of $A$ and coefficients $c(col_i)$; then to show that each summand $c(col_i) \langle\langle A,col_i\rangle\rangle$ is the weight of exactly one shadow-state. 

The  unknot $o$ is isotopic to $\partial N(T)$ and, while following the isotopy, at isolated moments it will cross some edges of $D\setminus (N(T)\cap D)$ (but no vertices or crossings because they are all contained in $N(T)$). Let us choose the isotopy so that every edge of $D\setminus (N(T)\cap D)$ is crossed exactly once (this can be done since each region is a disc and $T$ is a maximal connected sub-tree of $D$). 
We say that $u$ \emph{enters} a region $r$ if during the isotopy a subarc of $u$ not contained in $r$ crosses an edge contained in $r$. We claim that, since $T$, is connected each region $r_i,\ i=1,\ldots n$ will be ``entered" by $o$ exactly once during the isotopy. 
Indeed if $u$ enters twice a region $r_i$, let $\alpha,\beta$ the subarcs of $o$ in $r_i$, connecting them by an arc $\gamma$ we may produce two unknots whose connected sum is $u$. Since $T$ is connected and contains all the vertices and crossings of $D$, one of the two discs bounded by these unknots cannot contain vertices and so $\alpha$ and $\beta$ cross the same edge of $r_i$, against our hypothesis on the isotopy. 

\noindent\begin{minipage}{8cm}
We interpret each crossing moment as a fusion rule so that the isotopy of $u$ progressively ``erases" each arc of $D\setminus (N(T)\cap D)$ exactly when entering a region containing that arc, and the sum is taken over all admissible $u$: 
\end{minipage}
\begin{minipage}{7cm}
\begin{center}$$\raisebox{-.9cm}{\psfrag{a}{$a$}\psfrag{b}{$b$}\includegraphics[width=.9cm]{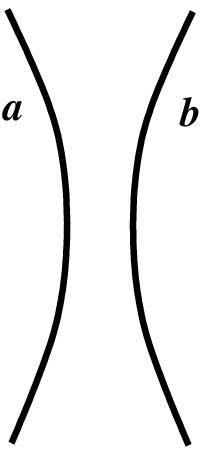}}={\Large{\sum}_{u}}\frac{\smallunknot{u}{.9}{.3}}{\smalltheta{a}{u}{b}{.9}{.3}}\raisebox{-.7cm}{\psfrag{a}{$a$}\psfrag{b}{$b$}\psfrag{c}{$u$}\includegraphics[width=.8cm]{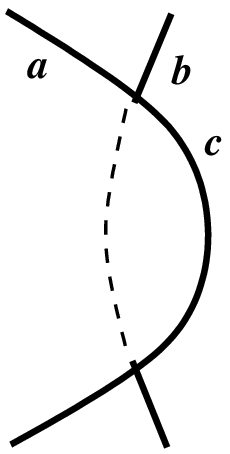}}$$
\end{center}
\end{minipage}

At the end of this isotopy, since for each $i\leq n$ $u$ entered $r_i$ only once, all the components of $\partial N(T)\cap r_i$ (which are arcs) will be colored by the same colors $u_i$ (see Figure \ref{fig:planargraph} for an example of this construction in the case of a planar graph). Therefore each summand in the final expression will be associated to a shadow-state $s$ given by $s(r_0)=0,s(r_i)=u_i$. Moreover, the other edges of $A$ (i.e. those of $N(T)\cap D$) are included in those of $G$ and therefore inherit the coloring $col$. Then we proved the following equality:
$$
\langle\langle G,col\rangle\rangle=\sum_{u_1,\ldots u_n} \prod_{r\in R} \smallunknot{u_i}{.9}{.4}\prod_{f\in F\setminus T}\smalltheta{}{}{}{.9}{.4}^{-1}\langle\langle A,col\cup \{u_1,\ldots u_n\}\rangle\rangle
$$
where the colors of the edges of the $\theta$ graphs are specified by the $u_i$'s and $col$. 
Remark that the summation range is exactly the set of shadow-states because the colors of the arcs of $\partial N(T)\cap r_i$ are all $u_i$ and the admissibility conditions for a three-uple of colors around an edge are satisfied at every moment we apply the fusion rule.
Moreover, in the above formula we already got part of the weights of each shadow-state (i.e. those of the regions and of all the edges out of $T$). 
We are left to prove that what is missing equals $\langle\langle A,col\cup \{u_1,\ldots u_n\}\rangle\rangle$, i.e. we claim the following:
$$ \langle\langle A,col(u_1,\ldots u_n)\rangle\rangle=\prod_{f\in T} \smalltheta{}{}{}{.9}{.4}^{-1}\prod_{v\in V} \smalltet{}{}{}{}{}{}{.9}{.4}\prod_{c\in C} \smalltetcross{}{}{}{}{}{}{.9}{.4}$$
where in each factor the colors are specified by the combinatorics of $D$ and by the state $u_1,\ldots u_n \cup col$ on $R\cup F$. To prove this, remark that: 
$$\langle\langle\raisebox{-.2cm}{\psfrag{a}{\footnotesize$a$}\psfrag{b}{\footnotesize$b$}\psfrag{c}{\footnotesize$u$}\includegraphics[width=1cm]{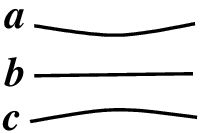}}\rangle\rangle=\sum_{i} \frac{\smallunknot{i}{.9}{.3}}{\smalltheta{a}{b}{i}{.9}{.3}}\ \langle\langle\raisebox{-.2cm}{\psfrag{a}{\footnotesize$a$}\psfrag{b}{\footnotesize$b$}\psfrag{c}{\footnotesize$u$}\includegraphics[width=1cm]{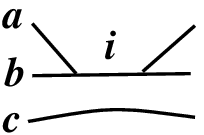}}\rangle\rangle=\smalltheta{a}{b}{c}{.9}{.3}^{-1}\ \langle\langle\raisebox{-.2cm}{\psfrag{a}{\footnotesize$a$}\psfrag{b}{\footnotesize$b$}\psfrag{c}{\footnotesize$u$}\includegraphics[width=1cm]{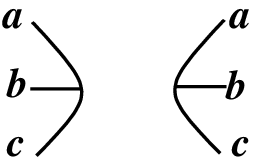}}\rangle\rangle$$
where the first equality is a fusion rule and the second is the inverse of a connected sum. Applying this identity on all the edges of $T$ we split $\langle\langle A,col\cup \{u_1,\ldots ,u_n\}\rangle\rangle$ into the product of the graphs remaining in the neighborhoods of the crossings and vertices (which are respectively \smalltetcross{}{}{}{}{}{}{.4}{.1} and \smalltet{}{}{}{}{}{}{.4}{.1}) divided by the product of the \smalltheta{}{}{}{.4}{.1}'s corresponding to the edges of $T$. This proves the claim and completes the proof when all the regions are discs. 
If in the beginning we added some $0$-colored edge to $G$ to cut $D$ into a diagram $D'$ whose regions are discs, then it is clear that every shadow state $s'$ of $D'$ can be lifted to a unique shadow state $s$ of $D$: indeed the compatibility conditions around a $0$-colored edge force the states of the neighboring regions to be the same.
Moreover, since the only difference between $D'$ and $D$ is given by the presence of the $0$-colored edge, it holds: 
$$w(s)=w(s')\prod_{f\in F_0} \smalltheta{u}{u}{0}{.9}{.4}^{-1}=w(s')\prod_{f\in F_0} \smallunknot{u}{.9}{.4}^{-1}=w(s')\prod_{r\in R} \smallunknot{u}{.9}{.4}^{\chi(r)-1}$$   
}\end{prf}

The formula given by Theorem \ref{teo:shadowstate} is often re-written by means of the so-called \emph{gleams}: 

\noindent\begin{minipage}{10cm}
\begin{defi}[Gleam]\label{def:gleam}
The gleam of a diagram $D$ of $G$ is a map $g:R\to \frac{\mathbb{Z}}{2}$ which on a region $r_i$ equals the sum over all the sectors of crossings contained in $r_i$ of $\frac{1}{2}$ times the local contributions of the crossing determined according to the pattern on the right. \end{defi}
\end{minipage}
\begin{minipage}{3cm}
\begin{center}
 \raisebox{0cm}{\psfrag{+}{$+$}\psfrag{-}{$-$}\includegraphics[width=1.5cm]{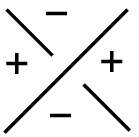}}
\end{center}
\end{minipage}
\begin{rem}
Do not confuse a shadow state with the gleam. In general for a given diagram $D$ there is a unique gleam but many different shadow-states.
\end{rem}
\begin{cor}\label{cor:finalshadowsum}
Under the same hypotheses as Theorem \ref{teo:shadowstate}, it holds:  
\begin{equation}\label{eq:shadowstatefinale}
\langle\langle G,col\rangle\rangle=F(G,col)\sum_{s}
\prod_{r\in R}  \sqrt{-1}^{4g(r)u}q^{2g(r)(u^2+u)}\smallunknot{u}{.9}{.4}^{\chi(r)}\prod_{f\in F} \smalltheta{}{}{}{.9}{.4}^{-\chi(f)}\prod_{V\cup C} \smalltet{}{}{}{}{}{}{.9}{.4} \end{equation}
\end{cor}
\begin{prf}{1}{
Using Example \ref{ex:crossedtet} one can rewrite the factors coming from crossings in term of tetrahedra multiplied by extra factors of the form $\sqrt{-1}^{\pm 2u}q^{\pm(u^2+u)}$ for each of the $4$ sectors around a crossing. To conclude, collect these factors according to the region containing the corresponding sectors and compare with Definition \ref{def:gleam}.
}\end{prf}
\subsection{Simplifying formulas}\label{sub:simpleformulas}
Both formulas \ref{eq:shadowstatefinale} and \ref{eq:shadowstate} are far from being optimal: indeed most of the factors in the state-sum can be discarded from the very beginning. Instead of giving a general theorem for doing this let us show why this happens through some examples.
We will say that an edge, vertex or crossing is \emph{external} if it is contained in the closure of $r_0$ and a region is external if its closure contains an external edge. 
Then the following simplifications can be operated:
\begin{enumerate}
\item If an external region contains two distinct external edges whose colors are different, then $\langle\langle G,col\rangle\rangle=0$ 
\item If an external region $r$ contains external edges $f_1,\ldots, f_k$ whose colors are all $c$, then, for every shadow-state $s$ on $D$, the total contribution coming from $r\cup f_1,\ldots f_k$ is: $$\smallunknot{c}{.9}{.3}^{\chi(r)-\sum_i \chi(f_i)}$$
\item If an external vertex $v$ is the endpoint of a non-external edge $f$ then their contribution simplify because: $$ \smalltet{a}{b}{c}{c}{0}{a}{.9}{.3}=\smalltheta{a}{b}{c}{.9}{.3}$$ (beware: if $\partial f$ is composed of two external vertices only one of them can be simplified with $f$)
\item For the same reason, if $D$ contains a sequence \raisebox{-.2cm}{\psfrag{a}{\footnotesize$a$}\psfrag{b}{\footnotesize$b$}\psfrag{r0}{\footnotesize$r_0$}\includegraphics[width=2cm]{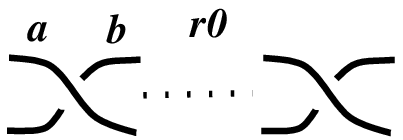}} of $n\in \mz\setminus \{0\}$ half-twists separating $r_0$ from a region $r$, then, in each summand of formula \ref{eq:shadowstatefinale}) the total contribution of the crossings and internal edges of the twist (excepted the initial and final ones) is $$\smalltheta{a}{b}{u}{.9}{.3}$$ where $u$ is the state of $r$ and $a,b$ are the colors of the strands (beware: the power of $q$ coming from the gleams do not simplify).
\end{enumerate}
\subsection{Examples and comments on integrality}\label{sub:commentsintegrality}
According to Theorem \ref{teo:main}, $\langle\langle G,col\rangle\rangle$ is a Laurent polynomial, and this is proved is by showing that the weight of each state in the state-sum expressing it via $R$-matrices and Clebsch-Gordan symbols is a Laurent polynomial.  Surprisingly enough, this is not true for shadow-state sums: the weight of a single shadow-state may be a rational function, but the poles of these functions will cancel out when summing on all the shadow-states. We will now clarify this by explict examples.
\subsubsection{Complicated formulas for unlinks}
Consider the $n$-colored unnormalized Jones polynomials of the following unlink:
$$
J_n(\raisebox{-.6cm}{\psfrag{+}{$\frac{1}{2}$}\psfrag{-}{-1}\includegraphics[width=3cm]{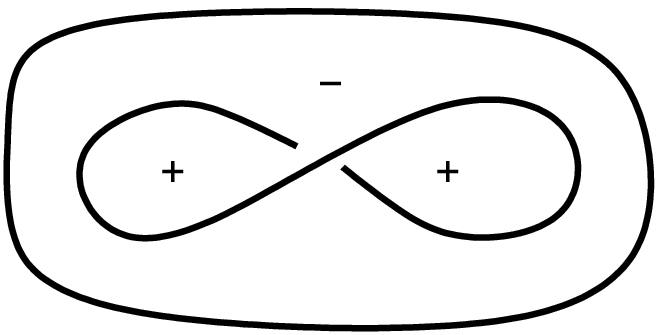}})= \sqrt{-1}^{-2n}q^{-4(n^2+n)}\sum_{u,v=0}^{2n} \sqrt{-1}^{2(u+v)}q^{(u^2+u)+(v^2+v)} \frac{\smallunknot{u}{.9}{.3}\smallunknot{v}{.9}{.3}\smalltet{n}{n}{u}{n}{n}{v}{.9}{.3}}{\smalltheta{n}{n}{u}{.9}{.3}\smalltheta{n}{n}{v}{.9}{.3}}
$$ 
In the picture we included the gleams of the regions to help the reader recovering formula \ref{eq:shadowstatefinale}). Of course this is a very complicated way of re-writing $(-1)^{2n}[2n+1]^2q^{2(n^2+n)}$, but what is interesting is that the single states are again not Laurent polynomials: for instance when $n=u=v=1$  the weight is $\frac{[3]([5]-1)}{[2][4]}$ .
\subsubsection{A more complicated link example}
Fix $a,b\in \mathbb{N}$ and consider the $n$-colored unnormalized Jones polynomials of the following link:
$$
J_n(\raisebox{-1cm}{\psfrag{n}{$a$}\psfrag{m}{$b$}\includegraphics[width=2cm]{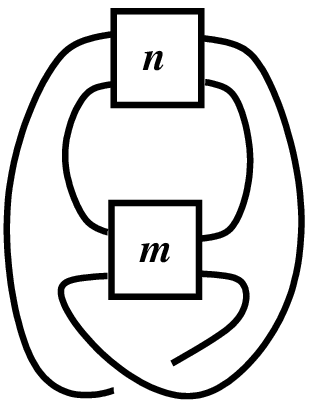}})=\sum_{u,v=0}^{2n}\sqrt{-1}^{2(a+1)u+2(b+1)v-4n(a+b+1)}q^{g(u,v)}\frac{\smallunknot{u}{.9}{.3}\smallunknot{v}{.9}{.3}\smalltet{n}{n}{u}{n}{n}{v}{.9}{.3}}{\smalltheta{n}{n}{u}{.9}{.3}\smalltheta{n}{n}{v}{.9}{.3}}$$
where a box with $a\in \mz$ stands for a sequence \raisebox{-0.1cm}{\psfrag{a}{}\psfrag{b}{}\psfrag{r0}{}\includegraphics[width=1.5cm]{twist.eps}} of $a$ half twists and $g(u,v)\doteqdot (a+1)(u^2+u)+(b+1)(v^2+v)-2(a+b+1)(n^2+n)$.
To get the above formula, before applying \ref{eq:shadowstatefinale}, remark that up to isotopy of the diagram, the region $r_0$ can chosen freely. Therefore, it is better to pick $r_0$ as the region touching the two boxes contemporaneously. 
Again, the summands are not Laurent polynomials but the sum is (take for instance $a=b=n=u=v=1$). More surprisingly, by Theorem \ref{teo:main}, for every $a,b\geq 0$ the resulting invariant will be divisible by $[2n+1]$ in $\mathbb{Z}[q,q^{-1}]$.   
\subsubsection{Some planar graphs}
If $G$ is a planar graph equipped with the blackboard framing, then the gleam of its regions are $0$ and so by Corollary \ref{cor:finalshadowsum} $\langle\langle G,col\rangle\rangle$ has a simple expression.
In the example of Figure \ref{fig:planargraph}, if all the edges of $G$ are colored by $n\in \mathbb{N}$ then
\begin{equation}\label{eq:cubone}
\langle\langle G,n\rangle\rangle=\sum_{u,v=0, \vert u-v\vert\leq n\leq u+v}^{2n} \smallunknot{u}{.9}{.3}\smallunknot{v}{.9}{.3}\frac{\smalltet{n}{n}{n}{u}{n}{n}{.9}{.3}^2\smalltet{n}{n}{n}{v}{n}{n}{.9}{.3}^2\smalltet{n}{n}{n}{u}{n}{v}{.9}{.3}^2}{\smalltheta{n}{n}{u}{.9}{.3}^3\smalltheta{n}{n}{v}{.9}{.3}^3\smalltheta{u}{n}{v}{.9}{.3}}
\end{equation}
where $u,v$ are the shadow-states of the two internal regions.
 If for instance in Formula \ref{eq:cubone} one puts $n=1$, and considers the shadow-state with $u=v=1$ then its weight is $\frac{[3]^2([4]!([5]-1))^6}{(-[4]!)^7}=-\frac{[3]^2([5]-1)^6}{[4]!}\notin \mathbb{Z}[q]$.
\subsection{Identities on $6j$-symbols}
Shadow state formulas provide a straightforward way to re-prove standard identities on $6j$-symbols. (The normalization we are using here for the symbols is that of Example \ref{ex:tet}).
\subsubsection{Normalizations of $6$-symbols}
It holds:
\begin{equation}
\delta_{b,0}\smallunknot{a}{.9}{.3}\smallunknot{c}{.9}{.3}=\sum_u \frac{\smallunknot{u}{.9}{.3}\smalltet{a}{a}{b}{c}{c}{u}{.9}{.3}}{\smalltheta{a}{c}{u}{.9}{.3}}
\end{equation}
where $u$ ranges between $\vert a-c\vert$ and $a+c$. This is proved by applying Formula \ref{eq:shadowstatefinale} to: \raisebox{-.5cm}{\psfrag{a}{\footnotesize$a$}\psfrag{b}{\footnotesize$b$}\psfrag{c}{\footnotesize$c$}\includegraphics[width=1.5cm]{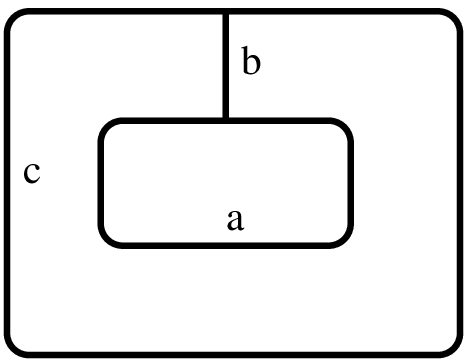}} and recalling that the invariant of a union of two unlinked graphs connected by a single arc is zero unless the color of the arc is $0$, in which case the invariant is just the product of the invariants of the graphs.
Similarly, it holds:
\begin{equation}
\smalltheta{a}{b}{c}{.9}{.3}\sqrt{-1}^{4b}q^{2(b^2+b)}=\sum_u \sqrt{-1}^{2(u+a-2c)}q^{u^2+u+a^2+a-2(c^2+c))}\frac{\smallunknot{u}{.9}{.3}\smalltet{a}{b}{c}{u}{b}{c}{.9}{.3}}{\smalltheta{b}{c}{u}{.9}{.3}}
\end{equation}
where $u$ ranges between $\vert b-c\vert$ and $b+c$. This is proved by applying Formula \ref{eq:shadowstatefinale} to \raisebox{-.5cm}{\psfrag{a}{\footnotesize$a$}\psfrag{b}{\footnotesize$b$}\psfrag{c}{\footnotesize$c$}\includegraphics[width=1.5cm]{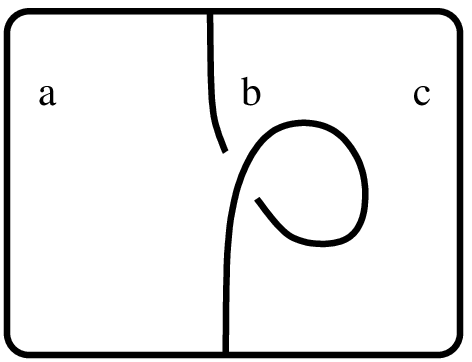}} and to the same \emph{framed} graph after undoing the kink on the $b$-colored edge.

\subsubsection{Orthogonality relation}
A direct corollary of Formula \ref{eq:shadowstatefinale} is the well-known orthogonality relation:
$$\delta_{b,d}\frac{\smalltheta{a}{b}{c}{.9}{.3}\smalltheta{e}{f}{b}{.9}{.3}}{\smallunknot{b}{.9}{.3}}=\sum_u \sqrt{-1}^{2(d-b)}q^{d^2+d-b^2-b}\frac{\smallunknot{u}{.9}{.3}\smalltet{a}{b}{c}{f}{u}{e}{.9}{.3}\smalltet{a}{c}{d}{f}{e}{u}{.9}{.3}}{\smalltheta{a}{e}{u}{.9}{.3}\smalltheta{c}{f}{u}{.9}{.3}}$$
where $u$ ranges over all the admissible colorings of the union of the two tetrahedral graphs on the right. To prove it, just apply formula \ref{eq:shadowstatefinale}) to the following two isotopic graphs and simplify the common factors:
$$\raisebox{-1cm}{\psfrag{a}{\footnotesize$a$}\psfrag{b}{\footnotesize$b$}\psfrag{c}{\footnotesize$c$}\psfrag{d}{\footnotesize$d$}\psfrag{e}{\footnotesize$e$}\psfrag{f}{\footnotesize$f$}\psfrag{g}{\footnotesize$g$}\psfrag{h}{\footnotesize$h$}\psfrag{i}{\footnotesize$i$}\includegraphics[width=3cm]{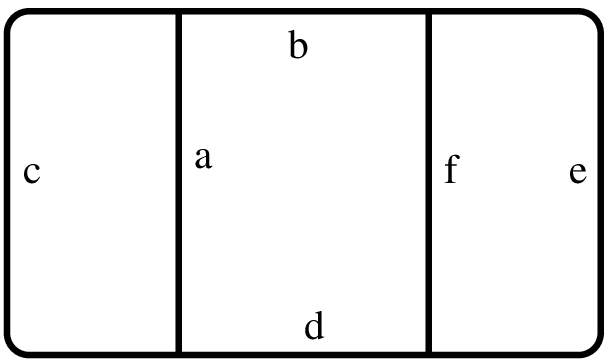}}\ \ \ \ \  \sim \ \ \ \ \ \ 
 \raisebox{-1cm}{\psfrag{a}{\footnotesize$a$}\psfrag{b}{\footnotesize$b$}\psfrag{c}{\footnotesize$c$}\psfrag{d}{\footnotesize$d$}\psfrag{e}{\footnotesize$e$}\psfrag{f}{\footnotesize$f$}\psfrag{g}{\footnotesize$g$}\psfrag{h}{\footnotesize$h$}\psfrag{i}{\footnotesize$i$}\includegraphics[width=3cm]{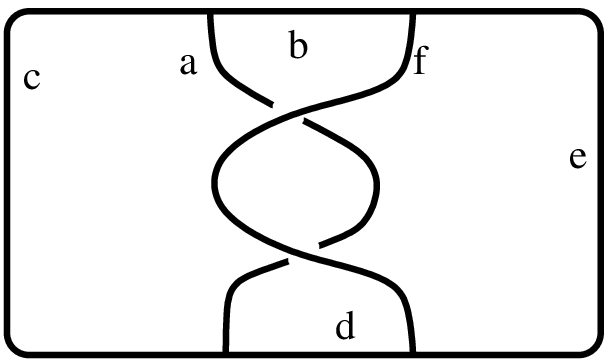}}$$

\subsubsection{Racah identity}
It holds:
\begin{equation}\label{eq:racah}
\sqrt{-1}^{2(a+b-c)}q^{a^2+a+b^2+b-c^2-c}\smalltet{a}{b}{c}{d}{e}{f}{.9}{.3}=\sum_{u}\sqrt{-1}^{2(u+f-e-d)}q^{u^2+u+f^2+f-d^2-d-g^2-e}\frac{\smallunknot{u}{.9}{.3}\smalltet{a}{e}{f}{b}{d}{u}{.9}{.3}\smalltet{a}{c}{b}{e}{u}{d}{.9}{.3}}{\smalltheta{a}{b}{u}{.9}{.3}\smalltheta{b}{d}{u}{.9}{.3}}
\end{equation}
where $u$ ranges over all the admissible colorings of the union of tetrahedra on the right.
Indeed it is sufficient to apply Formula \ref{eq:shadowstatefinale} to the following two isotopic graphs and simplify the common factors:
$$\raisebox{-1cm}{\psfrag{a}{\footnotesize$a$}\psfrag{b}{\footnotesize$b$}\psfrag{c}{\footnotesize$c$}\psfrag{d}{\footnotesize$d$}\psfrag{e}{\footnotesize$e$}\psfrag{f}{\footnotesize$f$}\psfrag{g}{\footnotesize$g$}\psfrag{h}{\footnotesize$h$}\psfrag{i}{\footnotesize$i$}\includegraphics[width=3cm]{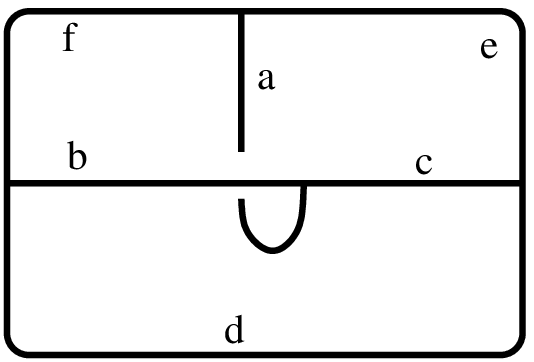}}\ \ \ \ \  \sim \ \ \ \ \ \ 
 \raisebox{-1cm}{\psfrag{a}{\footnotesize$a$}\psfrag{b}{\footnotesize$b$}\psfrag{c}{\footnotesize$c$}\psfrag{d}{\footnotesize$d$}\psfrag{e}{\footnotesize$e$}\psfrag{f}{\footnotesize$f$}\psfrag{g}{\footnotesize$g$}\psfrag{h}{\footnotesize$h$}\psfrag{i}{\footnotesize$i$}\includegraphics[width=3cm]{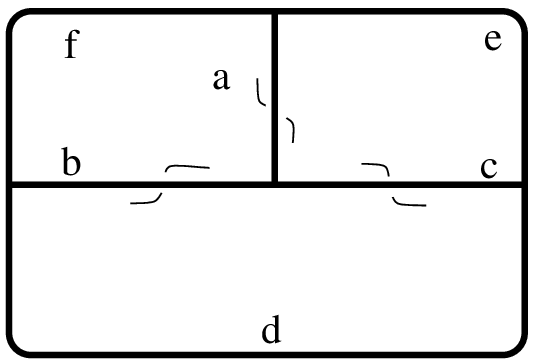}}$$

\subsubsection{Biedenharn-Elliot identity}
Another direct corollary of Formula \ref{eq:shadowstatefinale} is the well-known Biedenharn-Elliot identity:
\begin{equation}\label{eq:BE}
\frac{\smalltet{a}{b}{c}{d}{e}{f}{.9}{.3}\smalltet{g}{h}{e}{c}{d}{i}{.9}{.3}}{\smalltheta{c}{e}{d}{.9}{.3}}=\sum_{u}\frac{\smallunknot{u}{.9}{.3}\smalltet{a}{e}{f}{g}{u}{h}{.9}{.3}\smalltet{d}{b}{f}{u}{g}{i}{.9}{.3}\smalltet{a}{b}{c}{i}{h}{u}{.9}{.3}}{\smalltheta{a}{h}{u}{.9}{.3}\smalltheta{b}{i}{u}{.9}{.3}\smalltheta{f}{g}{u}{.9}{.3}}
\end{equation}
where $u$ ranges over all the admissible colorings of the union of tetrahedra on the right.
Indeed it is sufficient to apply Formula \ref{eq:shadowstatefinale} to the following two isotopic graphs and simplify the common factors:
$$\raisebox{-1cm}{\psfrag{a}{\footnotesize$a$}\psfrag{b}{\footnotesize$b$}\psfrag{c}{\footnotesize$c$}\psfrag{d}{\footnotesize$d$}\psfrag{e}{\footnotesize$e$}\psfrag{f}{\footnotesize$f$}\psfrag{g}{\footnotesize$g$}\psfrag{h}{\footnotesize$h$}\psfrag{i}{\footnotesize$i$}\includegraphics[width=3cm]{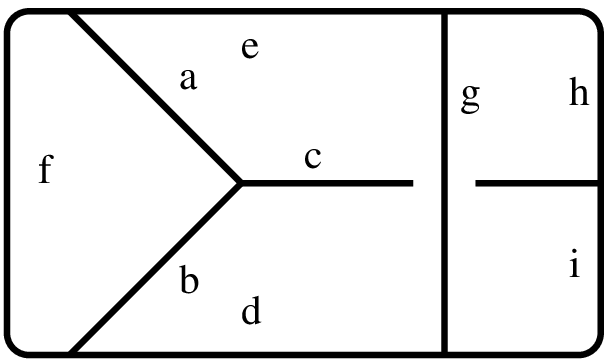}}\ \ \ \ \  \sim \ \ \ \ \ \ 
 \raisebox{-1cm}{\psfrag{a}{\footnotesize$a$}\psfrag{b}{\footnotesize$b$}\psfrag{c}{\footnotesize$c$}\psfrag{d}{\footnotesize$d$}\psfrag{e}{\footnotesize$e$}\psfrag{f}{\footnotesize$f$}\psfrag{g}{\footnotesize$g$}\psfrag{h}{\footnotesize$h$}\psfrag{i}{\footnotesize$i$}\includegraphics[width=3cm]{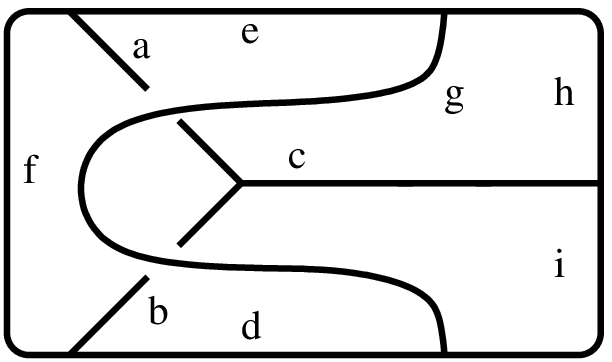}}$$

\section{$R$-matrices vs. $6j$-symbols}\label{sec:Rvs6j}
In the preceding sections we showed how to compute invariants of colored graphs by means of two different state-sums. Although for practical computation shadow state-sums turn out to be easier to deal with, state-sums based on $R$-matrices and Clebsch-Gordan symbols allow one to prove integrality results. 
In this section we will compare the state-sums when applied to graphs with boundary.
So let a $(n,m)$-KTG be a framed graph embedded in a square box which contains only $3$-valent vertices (inside the box) and $n$ (resp. $m$) $1$-valent vertices on the bottom (resp. top) edge of the box. A typical example is a framed $(n,m)$-tangle. Let as before $E,V$ be the set of edges and $3$-valent vertices of $G$, $\partial G^{\pm}\doteqdot \partial G\cap \partial box$ and, once chosen a diagram $D$ of $G$, let $F,C$ be the set of edges, and crossings of $D$. Let also $b_1,\ldots b_n$ be the bottom (univalent) vertices of $G$ and $t_1,\ldots t_m$ the top vertices. The definition of admissible coloring of a $(n,m)$ is the same as the standard one, but in this case, a second coloring is needed to get a numerical invariant out of $G$, namely a coloring on $\partial G$.
\begin{defi}[$\partial$-colorings for $(n,m)$-KTG's]
Let $(G,col)$ be a colored $(n,m)$-KTG; a $\partial$-coloring for $G$ is a map $col_{\partial}:\partial G\to \frac{\mathbb{Z}}{2}$ such that if $i_k$ (resp. $j_k$) is the color of the edge containing $b_k$ (resp. $t_k$) then $\vert col_{\partial}(b_k)\vert\leq i_k$ and $col_{\partial}(b_k)-i_k\in \mz$.
\end{defi}  
Equivalently a $\partial$-coloring is a choice of a vector in $V^{i_1}\otimes \cdots\otimes V^{i_n}$ of the form $g^{i_1}_{col_{\partial}(b_1)}\otimes \cdots \otimes g^{i_n}_{col_{\partial}(b_n)}$ and a vector in $V^{j_1}\otimes \cdots\otimes V^{j_m}$ of the form $g^{j_1}_{col_{\partial}(t_1)}\otimes \cdots \otimes g^{j_m}_{col_{\partial}(t_m)}$. 
Given a $(n,m)$-KTG equipped with a coloring $col$ and a $\partial$-coloring $col_{\partial}$, one can compute $\langle\langle G,col\cup col_{\partial} \rangle\rangle$ exactly as in Formula \ref{eq:statesum}: it is sufficient to restrict the set of admissible states to those such that the state of the boundary edges coincide with $col_{\partial}$.
Then $G$ represents a morphism $Z(G,col):V^{i_1}\otimes \cdots\otimes V^{i_n}\to V^{j_1}\otimes \cdots\otimes V^{j_m}$ and $\langle\langle G,col\cup col_{\partial}\rangle \rangle$ is an entry in the matrix expressing $Z(G,col)$ in the bases formed by tensor products of basis elements.

Most of the integrality result \ref{teo:main} still holds true (the idea of the proof is exactly the same):
\begin{teo}[Integrality: case with boundary]
The following belongs to $\mz[q^{\pm\frac{1}{2}}]$: $$\langle\langle G,col\cup col_{\partial}\rangle\rangle\doteqdot \langle G,col\cup col_{\partial}\rangle \frac{F(G,col)\prod_{e\in E'}[2col(e)]!}{\prod_{v\in V} [a_v+b_v-c_v]![a_v+c_v-b_v]![c_v+b_v-a_v]!}\frac{\prod_{k=1}^m\sqrt{-1}^{j_k}}{\prod_{k=1}^n\sqrt{-1}^{i_k}}$$ 
where $E'$ is the set of all the edges of $G$ which do not intersect $\partial G^+$ and $F(G,col)$ is defined as in the preceding sections.  
\end{teo}

What is interesting is that one can re-compute the invariant of $(G,col)$ also via shadow-state sums and Clebsch-Gordan symbols. To explain this, we will use the following:
\begin{defi}
Given a finite sequence $j_1,\ldots ,j_m$ a \emph{Bratteli sequence} associated to it is a sequence $s_0,s_1,\ldots s_m$ such that $s_0=0$ and for each $0\leq k\leq m-1$ the three-uple $s_k,j_k,s_k+1$ is admissible. 
\end{defi}
It is not difficult to realize that the set of Bratteli sequences associated to $j_1,\ldots ,j_m$ is in bijection with the set of irreducible submodules of $V^{j_1}\otimes \cdots\otimes V^{j_m}$ and that the submodule $V(s)$ corresponding to a Bratteli sequence $s=(s_0,\ldots s_m)$ is isomorphic to $V^{s_m}$. Moreover, using the morphisms defined in Subsections \ref{sub:Y} and \ref{sub:P}, we may fix explicit maps $\pi(s): V^{j_1}\otimes \cdots\otimes V^{j_m}\to V^{s_m}$ and $i(s):V^{s_m}\to V^{j_1}\otimes \cdots\otimes V^{j_m}$. Graphically they are expressed as follows: $$\pi(s)\doteqdot\raisebox{-.2cm}{\psfrag{s2}{\footnotesize $s_2$}\psfrag{sm}{\footnotesize $s_m$}\psfrag{j1}{\footnotesize $j_1$}\psfrag{j2}{\footnotesize $j_2$}\psfrag{j3}{\footnotesize $j_3$}\psfrag{jm}{\footnotesize $j_m$}\psfrag{dots}{\footnotesize $\cdots$}\includegraphics[width=2cm]{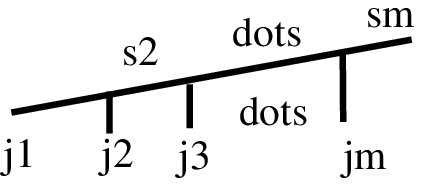}}\ \ \ \ \ \ \ \  {\rm and}\ \ \ \ \ \ \ i(s)\doteqdot\raisebox{-.2cm}{\psfrag{s2}{\footnotesize $s_2$}\psfrag{sn}{\footnotesize $s_m$}\psfrag{i1}{\footnotesize $j_1$}\psfrag{i2}{\footnotesize $j_2$}\psfrag{i3}{\footnotesize $j_3$}\psfrag{in}{\footnotesize $j_m$}\psfrag{dots}{\footnotesize $\cdots$}\includegraphics[width=2cm]{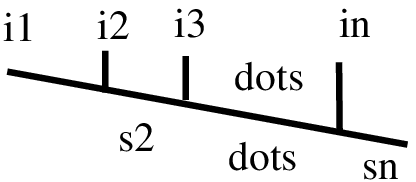}}$$

More explicitly, to construct $\pi(s)$, apply the Clebsch-Gordan morphism $P_{j_1,j_2}^{s_1}$ to the first two factors of $V^{j_1}\otimes \cdots \otimes V^{j_m}$ (the isomorphism is fixed by the choice of Clebsch-Gordan projectors made in Subsection \ref{sub:P}). The result is in $V^{s_2}\otimes V^{j_3}\otimes \cdots \otimes V^{j_m}$ and composing iteratively with with $P_{s_k,j_{k+1}}^{s_{k+1}}$ one gets the seeked map $\pi(s):V^{j_1}\otimes \cdots\otimes V^{j_m}\to V^{s_m}$. Similarly, $i(s):V^{s_m}\to V^{j_1}\otimes \cdots\otimes V^{j_m}$ is defined by composing recursively from right to left the Clebsch-Gordan morphisms of Subsection \ref{sub:Y} $V^{s_m}\to V^{s_{m-1}}\otimes V^{j_m}\to V^{s_{m-2}}\otimes V^{j_{m-2}}\otimes V^{j_m}\to\ldots \to V^{j_1}\otimes \cdots\otimes V^{j_m}$.
  
\begin{figure}
\psfrag{s+1}{$s^+_1$}
\psfrag{s-1}{$s^-_1$}
\psfrag{s+m}{$s^+_{m-1}$}
\psfrag{s-n}{$s^-_{n-1}$}
\psfrag{x}{$x$}
\psfrag{0}{$0$}
\psfrag{G}{$G$}
\psfrag{dots}{$\cdots$}
\includegraphics[width=13cm]{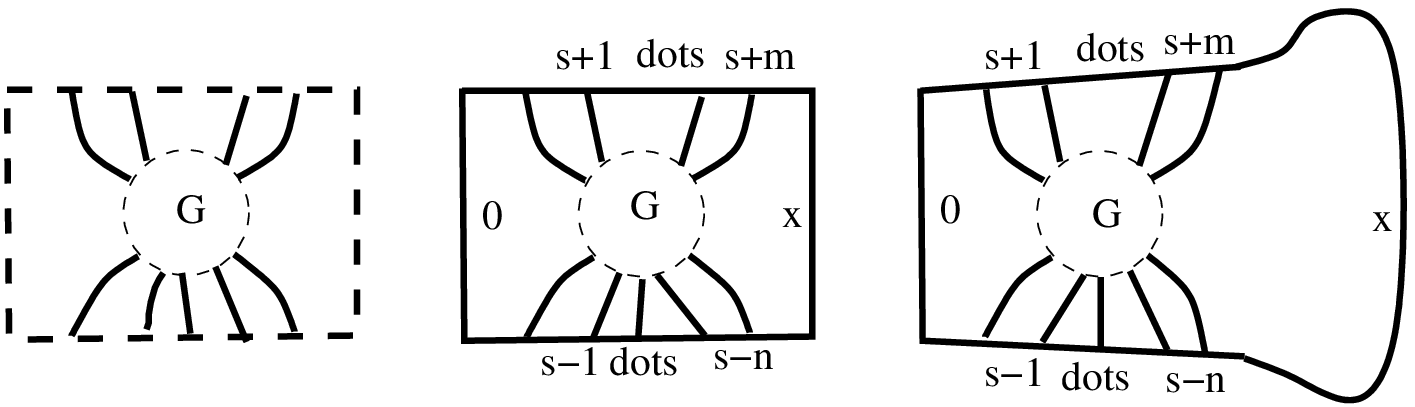}
\caption{Using a pair of Bratteli sequences it is possible to close $(G,col)$ to a colored closed $KTG$ (in the middle).}\label{fig:bratteli}
\end{figure}

Let  $\overline{G}\doteqdot G\cup \partial box$ viewed as a framed graph in $\mr^2$ by embedding the box in $\mr^2$ with the blackboard framing around its boundary. 
Let $s^+=(s^+_1,\ldots s^+_m)$ and $s^-=(s^-_1,\ldots s^-_m)$ be Brattelli sequences associated to $j_1,\ldots, j_m$ and $i_1,\ldots ,i_n$ and suppose that $s^+_m=s^-_n\doteqdot x$.
We can extend $col$ to a coloring $col\cup s^-\cup s^+$ of $\overline{G}$: the color of the edge in the top (resp. bottom) edge of the box bounded by $j_k$ and $j_{k+1}$ (resp.  $i_k$ and $i_{k+1}$) is $s^+_k$ (resp. $s^-_k$), the color of the left edge of the box is $0$ and that of the right edge $x$ (see Figure \ref{fig:bratteli}).
\begin{teo} [Shadow state-sums vs. R-matrices and Clebsch-Gordan symbols]
Let $\lambda(s^-,s^+)\in \mc$ be defined by $\pi(s^+)\circ Z(G,col)\circ i(s^-)=\lambda(s^-,s^+) Id_{V^x}$. Then:
$$\lambda(s^-,s^+)\smallunknot{x}{.9}{.3}=\langle\langle \overline{G},col\cup s^-\cup s^+\rangle\rangle$$
and
$$Z(G,col)=\sum_{s^-}\sum_{s^+} \prod_{t=1}^{n-1} \frac{\smallunknot{s_t}{.9}{.3}}{\smalltheta{i_t}{s_{t}}{\nns \nns \raisebox{.1cm}{$s_{t+1}$}}{.9}{.3}}  \prod_{t=1}^{m-1}\frac{\smallunknot{s_t}{.9}{.3}}{\smalltheta{j_t}{s_{t}}{\nns \nns \raisebox{.1cm}{$s_{t+1}$}}{.9}{.3}} \frac{\langle\langle \overline{G},col\cup s^-\cup s^+\rangle\rangle}{\smallunknot{s_n}{.9}{.3}}\raisebox{-1.0cm}{\psfrag{s2}{\footnotesize$s^+_2$}\psfrag{s-2}{\footnotesize$s^-_2$}\psfrag{s-n}{\footnotesize$s^-_n$}\psfrag{sm}{\nns\footnotesize$s^+_m$}\psfrag{i1}{\footnotesize$i_1$}\psfrag{i2}{\footnotesize$i_2$}\psfrag{i3}{\footnotesize$i_3$}\psfrag{in}{\footnotesize$i_n$}\psfrag{j1}{\footnotesize$j_1$}\psfrag{j2}{\footnotesize$j_2$}\psfrag{j3}{\footnotesize$j_3$}\psfrag{jm}{\footnotesize$j_m$}\psfrag{dots}{\footnotesize$\cdots$}\includegraphics[width=3cm]{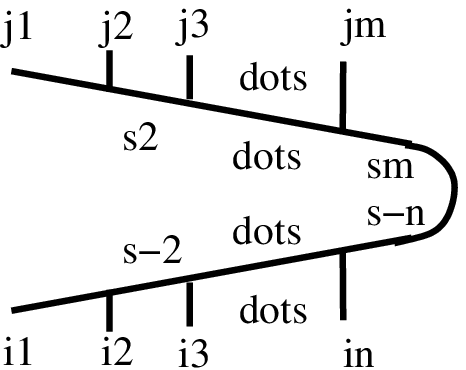}}.
$$
where the sum are taken over all the Bratteli sequences $s^-$ and $s^+$ associated respectively to $i_1,\ldots i_n$ and $j_1,\ldots j_m$.
\end{teo}

\begin{prf}{1}{
The value on the right hand side is the invariant of the colored graph $(\overline{G},col\cup s^+\cup s^-)$ depicted on the right of Figure \ref{fig:bratteli}. But, as shown in the picture, $\overline{G}$ is the closure of the graph representing $\pi(s^+)\circ Z(G,col)\circ i(s^-)$.

The second statement follows by applying $n-1$ times the fusion rule on the bottom legs of $G$ and $m-1$ times on the top legs to get:\\
\begin{center}
\raisebox{-1.5cm}{\psfrag{s+1}{$s^+_1$}\psfrag{s-1}{$s^-_1$}\psfrag{s+m}{$s^+_{m-1}$}\psfrag{s-n}{$s^-_{n-1}$}\psfrag{x}{$x$}\psfrag{0}{$0$}\psfrag{G}{$G$}\psfrag{dots}{$\cdots$}\includegraphics[width=3cm]{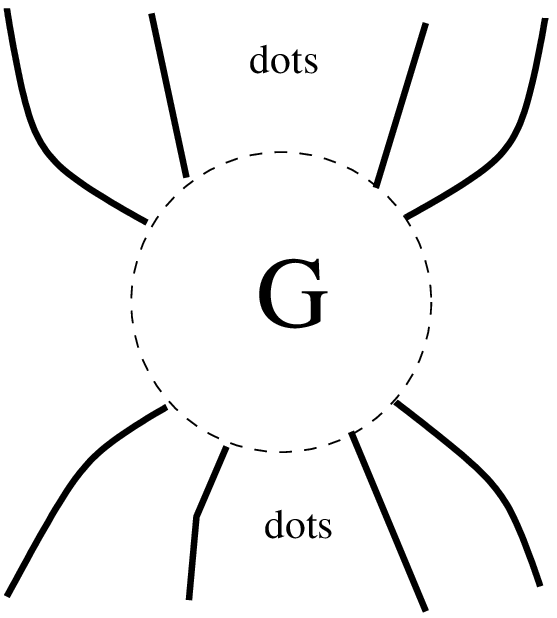}}=$
\sum_{s^-}\sum_{s^+} \prod_{t=1}^{n-1} \frac{\smallunknot{s_t}{.9}{.3}}{\smalltheta{i_t}{s_{t}}{\nns \nns \raisebox{.1cm}{$s_{t+1}$}}{.9}{.3}}  \prod_{t=1}^{m-1}\frac{\smallunknot{s_t}{.9}{.3}}{\smalltheta{j_t}{s_{t}}{\nns \nns \raisebox{.1cm}{$s_{t+1}$}}{.9}{.3}}$
\raisebox{-1.5cm}{\psfrag{s2}{\footnotesize$s^+_2$}\psfrag{s-2}{\footnotesize$s^-_2$}\psfrag{s-n}{\footnotesize$s^-_n$}\psfrag{sm}{\nns \nns \footnotesize$s^+_m$}\psfrag{i1}{\footnotesize$i_1$}\psfrag{i2}{\footnotesize$i_2$}\psfrag{i3}{\footnotesize$i_3$}\psfrag{in}{\footnotesize$i_n$}\psfrag{j1}{\footnotesize$j_1$}\psfrag{j2}{\footnotesize$j_2$}\psfrag{j3}{\footnotesize$j_3$}\psfrag{jm}{\footnotesize$j_m$}\psfrag{dots}{\footnotesize$\cdots$}\psfrag{GB}{$\overline{G}'$}\psfrag{dots}{$\cdots$}\includegraphics[width=3cm]{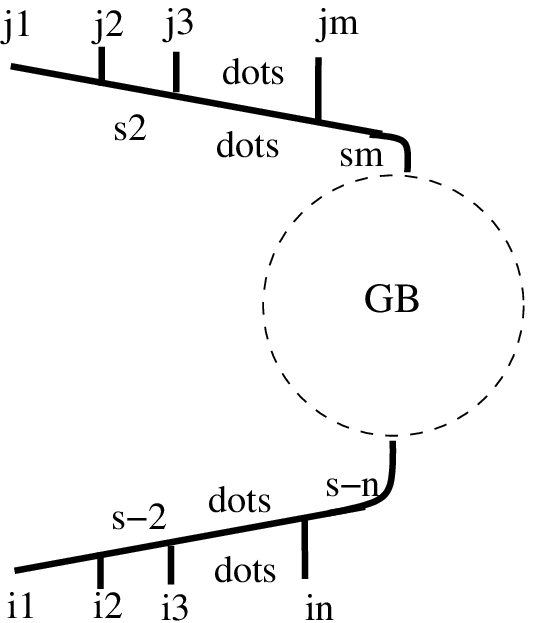}}
\end{center}
where $s^-$ and $s^+$ range over all Bratteli sequences associated to $i_1,\ldots ,i_n$ and $j_1,\ldots ,j_n$ respectively and $\overline{G}'$ is the $(1,1)$ colored KTG obtained by opening  $(\overline{G}, col \cup s^+\cup s^-)$ along the right edge of the box.
}\end{prf}
\begin{example}[$R$-matrices vs $6j$-symbols]
It holds: $${}_b^aR^{t,w}_{u,v}=\sum_{c=\vert a-b\vert}^{a+b} C^{b,a,c}_{t,w,v+u}\frac{\sqrt{-1}^{2(c-a-b)}q^{c^2+c-a^2-a-b^2-b}\smallunknot{c}{.9}{.3}}{\smalltheta{a}{b}{c}{.9}{.3}}P^{a,b,c}_{u,v,u+v}$$
Indeed it is sufficient to apply the preceding theorem to $(G,col)$ being a crossing. 
\end{example}

\end{document}